\documentclass[12pt]{article}
\usepackage{amsmath,amsfonts,amssymb,amsthm}
\usepackage{hyperref}
\usepackage{graphicx}

\newcommand{\ind}{{\bf 1}}
\newtheorem{proposition}{Proposition}[section]
\newtheorem{theorem}[proposition]{Theorem}
\newtheorem{corollary}[proposition]{Corollary}
\newtheorem{lemma}[proposition]{Lemma}
\newtheorem{remark}[proposition]{Remark}

\numberwithin{equation}{section}
\newcommand{\nc}{\newcommand}
\nc{\R}{{\mathbb R}}
\newcommand{\I}{{\bf 1}}

\nc{\G}{{\mathbb G}}
\nc{\cG}{{\mathcal G}}
\nc{\N}{{\mathbb N}}
\nc{\Z}{{\mathbb Z}}
\renewcommand{\P}{\mathbb{P}}
\nc{\E}{\mathbb{E}}
\nc{\Q}{\mathbb{Q}}
\renewcommand{\H}{\mathbb{H}}
\nc{\bN}{{\mathbf N}}
\nc{\bM}{{\mathbf M}}
\nc{\X}{{\mathbb X}}
\nc{\Y}{{\mathbb Y}}
\nc{\bx}{{\mathbf x}}
\nc{\cB}{{\mathcal B}}
\nc{\cF}{{\mathcal F}}

\DeclareMathOperator*{\sumneq}{\sideset{}{^{\neq}}\sum}

\nc{\inclD}[1]{\vcenter{\hbox{\includegraphics{Images/#1}}}}

\begin{document}

\author{G\"unter Last\thanks{Karlsruhe Institute of Technology, guenter.last@kit.edu}
 and Sebastian Ziesche\thanks{Karlsruhe Institute of Technology, sebastian.ziesche@kit.edu}
}

\title{On the Ornstein-Zernike equation for\\stationary cluster processes and the\\
random connection model}
\date{\today}
\maketitle

\begin{abstract}
\noindent In the first part of this paper we consider
a general stationary subcritical cluster model in $\R^d$.
The associated pair-connectedness function can be defined
in terms of two-point Palm probabilities of the underlying point
process. Using Palm calculus and Fourier theory we
solve the Ornstein-Zernike equation (OZE) under quite general
distributional assumptions.
In the second part of the paper we discuss the
analytic and combinatorial properties of the OZE-solution in the special
case of a Poisson driven random  connection model.
\end{abstract}

\maketitle

\begin{flushleft}
\textbf{Key words:} Ornstein-Zernike equation, Gilbert graph,
random connection model, pair-connectedness function, Poisson process, percolation, analyticity
\newline
\textbf{MSC (2010):} Primary: 60K35, 60D05; Secondary: 60G55
\end{flushleft}

\section{Introduction}\label{secintro}

In a seminal paper Ornstein and Zernike
proposed in \cite{OZ14} to split the interaction
between molecules in a liquid into a direct
and an indirect part.
While the resulting spatial convolution equation is
of great important in physics, it seems to be
hardly known among mathematicians. The
aim of this paper is to bridge this gap and to
lay a rigorous mathematical foundation for further studies.

We start with a simple example of a stationary cluster process, which is also
a special case of the random connection model, studied later.
Let $\eta_t$ be a stationary Poisson process on $\R^d$ with intensity $t\ge 0$.
Let $B\subset\R^d$ be a {\em gauge body}, that is a
compact 
set containing the origin $0\in\R^d$ in its interior.  We define a
{\em random geometric graph} $\cG(\eta_t)$ with vertex set $\eta_t$ as follows.
Two distinct points
$x,y\in\eta_t$ are adjacent in $\cG(\eta_t)$ whenever $(B+x)\cap (B+y)\ne\emptyset$,
where $B+x:=\{x+z:z\in B\}$; see \cite{Pen03}.
For $x\in\eta_t$ let $C(x,\eta_t)\subset \eta_t$ denote the
{\em cluster} of $x$, that is the connected
component of $\cG(\eta_t)$ containing $x$.
These definitions apply to any point process $\eta$,
and in particular to deterministic locally finite subsets of $\R^d$.
For each point process $\eta$ on $\R^d$ and $x,y\in\R^d$ we write
$\eta^x:=\eta\cup\{x\}$ and $\eta^{x,y}:=\eta\cup\{x,y\}$; see also the Appendix.
We wish to study the {\em pair connectedness function} (see \cite{Tor02})
\begin{align}\label{pairconn}
P_t(x,y):=\P(y\in C(x,\eta^{x,y}_t)),\quad x,y\in\R^d.
\end{align}
By Corollary~4.15 in \cite{Hall88}, there is a
{\em percolation threshold} $t_c\in(0,\infty)$ such that
$\P(|C(0,\eta^0_t)|=\infty)>0$ for $t>t_c$ and
$\P(|C(0,\eta^0_t)|=\infty)=0$ for $t<t_c$.
We seek a function $Q_t(x,y)$ solving the {\em Ornstein-Zernike equation}
\begin{align}\label{OZ1}
  P_t(x,y)=Q_t(x,y)+t\int_{\R^d} Q_t(x,z)P_t(z,y)\,dz,\quad x,y\in\R^d,\,t<t_c.
\end{align}
We shall formulate and solve equation \eqref{OZ1} in the following
much more general setting. Let $\eta$ be a stationary point
process on $\R^d$ with finite intensity $\gamma_\eta$. The points are
partitioned into clusters (sets of points of $\eta$) according to a
translation-invariant rule. This rule might be very general and can
incorporate additional randomness (e.g.\ in the random connection
model).
The point process $\eta$ is assumed to be jointly
stationary with the cluster process.  The pair-connectedness function
$P(x,y)$ is then informally defined as the conditional probability
that $x,y\in\R^d$ belong to the same cluster, given that $x$ and $y$
are points of $\eta$ (suitably weighted by the pair-correlation
function). Then the Ornstein-Zernike equation \eqref{OZ1} takes the
form
\begin{align}\label{OZ3}
  P(x,y)=Q(x,y)+\gamma_\eta\int_{\R^d} Q(x,z)P(z,y)\,dz,\quad x,y\in\R^d.
\end{align}
Our Theorem \ref{tOZ} shows under rather weak assumptions
that equation \eqref{OZ3} has a unique solution. The proof of this result
is based on Palm calculus
for stationary point processes (see the Appendix) and
a classical theorem by Wiener on the inversion of Fourier transforms.

In Sections \ref{secrcm}--\ref{seccomb} of the paper we shall consider
the (Poisson driven) random connection model (RCM) (see
\cite{MeesterRoy96}), a significant generalization of the Gilbert
graph introduced above. The RCM with parameters $t \geq 0$ and
$\varphi\colon\R^d\to[0,1]$ is a random graph $\cG$ where the set of
vertices is a Poisson process $\eta_t$ with intensity $t$. Any two
distinct vertices $x,y \in \eta_t$ are adjacent with probability
$\varphi(x-y)$ independently of all other pairs and of $\eta_t$. We call
$\varphi$ the \emph{connection function} of the RCM.
The clusters in this model are just the connected components
of  $\cG$. In Section
\ref{secrcm} we shall give a detailed description of this model along
with formulas on degree distributions (that are basically well-known)
and a Margulis-Russo type formula. The latter result might be of some
independent interest. In Section \ref{secOZErcm} we shall first show
that 
the RCM satisfies the assumptions of Theorem \ref{tOZ},
so that a solution $Q_t\equiv Q$ of \eqref{OZ1} (with $P_t\equiv P$ denoting the pair
connectedness function) exists in the whole subcritical regime.  Then
we prove that $P_t$ is an analytic function of $t$ on the interval
$[0,t_*)$, where $t_*$ is the smallest number such that for $t<t_*$
the typical cluster has an exponentially decreasing tail.  In the
Boolean model with fixed gauge body $B$ (mentioned above), the
arguments from \cite{Pen03} can be extended to show that $t_*$ is
equal to the percolation threshold $t_c$.
In fact, Theorem \ref{analyticityOfClusterFctls} shows that this
result holds for general integrable functions of the typical cluster.
We then proceed with deriving similar properties for $Q_t$; see
Proposition \ref{analyticity-Q_t}.  We are not aware of a direct
probabilistic interpretation of
$Q_t$. However, for small intensities $t$ there
is a simple combinatorial relationship between
the coefficients in the expansions of $P_t$ and $Q_t$;
see Theorem \ref{thm:repr_q_n}.

In writing this paper we strongly benefited from the large physics literature
on the topic.
In particular the combinatorial formulas in the final section
are basically well-known, although not in a mathematically rigorous form.
Two key references are \cite{Bax70} and \cite{Con77}.
However, we have not been able to find a justification
of the existence of solution of the OZE, even not
in the very special case of the Poisson driven Boolean model.
Moreover, the analytic properties of $P_t$ and $Q_t$ (often taken
as granted) have not been proved either.
In our opinion it is one of the main contributions of the present paper to
apply modern point process methods (Palm calculus and
Margulis-Russo type formulas for Poisson driven systems) to the
OZE.
The original motivation for our work came from \cite{Tor12}, where the
author uses the OZE to derive putative lower bounds
for the percolation threshold in the Boolean model. This is one of the potential
applications of the present paper.

\section{The Ornstein-Zernike equation}\label{secOZ}

In this section we establish equation \eqref{OZ3} for general
stationary cluster processes defined on a probability space $(\Omega,\mathcal{A},\P)$.
As in the Appendix we assume that $(\Omega,\mathcal{A})$ is equipped
with a measurable flow $\{\theta_x:x\in\R^d\}$ leaving $\P$
invariant.
We let $\eta$ be an invariant (and therefore stationary) point process on $\R^d$
with finite intensity $\gamma_\eta:=\E\eta([0,1]^d)$.
We also assume that $\P(\eta\ne\emptyset)=1$.
To describe the clusters, we consider a measurable mapping
$(\omega,x)\mapsto \tau(\omega,x)$ from $\Omega\times\R^d$ into
$\R^d$ with the covariance property
\begin{align}\label{covariance}
\tau(\theta_y\omega,x-y)=\tau(\omega,x)-y,\quad \omega\in\Omega,\, x,y\in\R^d.
\end{align}
(For convenience we also assume that $\tau(x)=x$, $x\in\R^d$, whenever $\eta(\R^d)=0$.)
The points of the random set
\begin{align*}
\xi:=\{\tau(x):x\in\eta\}
\end{align*}
are interpreted as locations (or centers) of the clusters of $\eta$.
Note that $\tau$ need not to be a
deterministic function of $\eta$ like in the Boolean model, but might
incorporate additional randomness; see Section \ref{secrcm}.
The refined Campbell formula \eqref{refC} and the
covariance propert \eqref{covariance} imply that
\begin{align*}
  \E|\xi\cap B|=\gamma_\eta\lambda_d(B)\,\P^0_\eta(\tau(0)\in B),\quad B\in\cB^d,
\end{align*}
where $|A|$ denotes cardinality of a set $A$,
$\lambda_d$ is Lebesgue measure on $\R^d$ and $\P^0_\eta$
is the Palm probability measure of $\eta$; see Appendix.
In particular $\E|\xi\cap B|<\infty$ for all bounded Borel sets $B$,
so that it is no restriction of generality to assume that
$\xi$ is locally finite everywhere on $\Omega$. It follows that
$\xi$ is an invariant point process with finite intensity
$\gamma_\xi=\gamma_\eta\,\P^0_\eta(\tau(0)\in [0,1]^d)$.
The {\em clusters} can be formally defined as those
points of $\eta$ which have the same image under $\tau$.
Hence $x,y\in\eta$ belong to the same cluster
iff $\tau(x)=\tau(y)$ and the cluster of $x\in\eta(\omega)$ is given by
$C(\omega,x):=\{y\in \eta(\omega):\tau(\omega,y)=\tau(\omega,x)\}$ or,
more succinctly,
\begin{align}\label{e2.3}
C(x):=\{y\in \eta:\tau(y)=\tau(x)\},\quad x\in\eta.
\end{align}
(It is convenient to use this definition for all $x\in\R^d$.)
In the random connection model, for instance, the mapping $\tau$
is defined so as to make sure that \eqref{e2.3}
is consistent with the definition of the clusters given in
the introduction.

It follows from \eqref{covariance} that
\begin{align}\label{e2.4}
C(\theta_y\omega,x)=C(\omega,x+y)-y,\quad \omega\in\Omega,\, x,y\in\R^d.
\end{align}
The distribution of $C(0)$ under the Palm probability measure $\P^0_\eta$ (see Appendix)
can be interpreted as the distribution of the cluster containing
the typical point of $\eta$. We make the crucial assumption that the size
of this cluster has a finite expectation, that is
\begin{align}\label{e2.5}
\E^0_\eta |C(0)|<\infty,
\end{align}
where $\E^0_\eta$ denotes the expectation with respect to $\P^0_\eta$.
To retrieve the points in a cluster with location $z \in \xi(\omega)$ we define
$D(\omega,z):=\{x\in\eta(\omega):\tau(\omega,x)=z\}$ or,
\begin{align}\label{e2.6}
D(z):=\{x\in \eta:\tau(x)=z\},\quad z\in\xi.
\end{align}
(Again we use this notation for all $z\in\R^d$.)

As we are interested in second order properties of $\eta$, we need
to assume that the second order factorial moment measure of $\eta$ is
locally finite and absolutely continuous. We then denote by $\rho$ the pair correlation
function and by $\P^{x,y}_\eta$, $x,y\in\R^d$, the bivariate Palm distributions
of $\eta$. The latter are probability measures on $(\Omega,\mathcal{A})$;
see Appendix. We can interpret $\P^{x,y}_\eta(A)$
as the conditional probability of $A\in\mathcal{A}$ given that $\eta$ has points
at $x$ and $y$. Our interest in this paper focuses
on the weighted {\em pair connectedness function}
\begin{align}\label{e2.10}
P(x,y):=\rho(x-y)\P^{x,y}_\eta(y\in C(x)),\quad x,y\in\R^d.
\end{align}
In view of  \eqref{A3} and \eqref{e2.4} we have $P(x,y)=P(y-x)$ and
we define the (even) function $P\colon \R^d\to\R$ by
$P(x):=P(0,x)$. Choosing $f:=\I\{0\in C(x)\}$ in \eqref{A5} gives
\begin{align}\label{e2-12}
\E^0_\eta|C(0)|=1+\gamma_\eta\int P(x)\,dx.
\end{align}
Hence \eqref{e2.5} implies that $P$ is in the space
$L^1$ of all measurable functions $f\colon \R^d\to\R$ with
$\|f\|_1:=\int |f(x)|\,dx<\infty$.

The convolution of $f,g\in L^1$ is defined
as
$$
(f\ast g)(x):=\int f(x-y)g(y)\,dy,\quad x\in\R^d.
$$
In the same way we define the convolution for functions $f \in L^1$
and $g \in L^\infty $ where $L^\infty $ is the space of bounded
functions equipped with the supremum norm $\|\cdot\|_\infty$. Both definitions make sense due to
the basic inequalities
\begin{align}
  &\|f \ast g\|_\infty \leq \|f\|_1 \|g\|_\infty , \quad f \in L^1, g \in L^\infty,\label{young-ineq-inf}\\
  &\|f \ast g\|_1 \leq \|f\|_1 \|g\|_1 , \quad f \in L^1, g \in L^\infty \label{young-ineq-1}.
\end{align}

We can now formulate and prove the Ornstein-Zernike equation
\eqref{OZ1} in the present very general stationary setting.
We need the  regularity assumption
\begin{align}\label{nozeroeta}
\P^0_{\eta}\bigg(\sum_{x\in C(0)} e^{\mathbf{i}wx}\ne 0\bigg)>0,\quad w\in\R^d,
\end{align}
where $wx$ is the Euclidean scalar product of $x,w\in\R^d$ and
$\mathbf{i}$ is the imaginary unit.

\begin{theorem}\label{tOZ} Assume that \eqref{e2.5} and \eqref{nozeroeta} hold.
Then there is a unique $Q\in L^1 \cap L^\infty $ such that
\begin{align}\label{OZ}
P=Q+\gamma_\eta Q\ast P.
\end{align}
\end{theorem}

\begin{remark}\label{rassumption}\rm Assumption \eqref{nozeroeta} is rather weak.
It holds, for instance, if $\P^0_{\eta}(|C(0)|=1)>0$. Indeed, if
$|C(0)|=1$, then the sum in \eqref{nozeroeta} reduces to the single term $1$.
Another sufficient condition can be formulated in terms of the
factorial moment measures $\alpha^{(n)}$, $n\in\N$, of $\eta$  defined by
\eqref{factorial}.
If these measures are locally finite and absolutely continuous, then
\begin{align*}
\P^0_{\eta}\bigg(\sum_{x\in C(0)} e^{\mathbf{i}wx}= 0\bigg)=0,\quad w\in\R^d,
\end{align*}
so that \eqref{nozeroeta} holds. To see this, we note that
\begin{align*}
\P^0_{\eta}\bigg(\sum_{x\in C(0)} e^{\mathbf{i}wx}= 0\bigg)
&\le \E \sum^\infty_{n=2}\frac{1}{n!}\sumneq_{x_1,\ldots,x_n\in\eta}
\I\{e^{-\mathbf{i}wx_1}f_n(x_2,\ldots,x_n)=1\}\\
&=\sum^\infty_{n=2}\frac{1}{n!}\E \int\I\{f_n(x_2,\ldots,x_n)= e^{\mathbf{i}wx_1}\}\,\alpha^{(n)}(d(x_1,\ldots,x_n)),
\end{align*}
where $f_n(x_2,\ldots,x_n):=\sum^n_{k=2} e^{\mathbf{i}wx_k}$.

\end{remark}

We prepare the proof of Theorem \ref{tOZ} with some
results of independent interest. We start with the classical
connection between the typical cluster and the cluster of a typical point.
In what follows we interpret $C(x)$ and $D(x)$ as point processes
on $\R^d$, i.e.\ as measurable mappings from $\Omega$ to $\mathbf{N}(\R^d)$.

\begin{proposition}\label{prop2.1} For any measurable $f\colon \bN(\R^d)\to[0,\infty)$
\begin{align}\label{e2.15}
\gamma_\eta\, \E^0_\eta f(C(0)-\tau(0))=\gamma_\xi \, \E^0_{\xi} |D(0)|f(D(0))
\end{align}
\end{proposition}
\emph{Proof:} We have that
\begin{align*}
\gamma_\xi \, \E^0_{\xi} |D(0)|f(D(0))
&=\gamma_\xi \, \E^0_{\xi}\sum_{x\in\eta} f(D(0))\I\{\tau(x)=0\}\\
&= \gamma_\eta\, \E^0_{\eta}\sum_{x\in\xi} f(D(\theta_x,0))\I\{\tau(\theta_x,-x)=0\}\\
&= \gamma_\eta\, \E^0_{\eta}\sum_{x\in\xi} f(D(\theta_x,0))\I\{\tau(\theta_0,0)=x\},
\end{align*}
where we have used Proposition \ref{tneveu} with
$(\omega,x) \mapsto f(D(\omega, 0))\I\{\tau(\omega,x)=0\}$ to get the second and
\eqref{covariance} to get the third identity. Using \eqref{covariance} again,
it can be easily checked that $D(\theta_x\omega,0)=C(\omega,0)-x$,
whenever $x\in\xi(\omega)$ and $\tau(\omega, 0) = x$. This finishes the proof.\qed

\bigskip
Proposition \eqref{prop2.1} implies in particular that
\begin{align}\label{e2.12}
 \gamma_\xi\, \E^0_{\xi}|D(0)|^2=\gamma_\eta\, \E^0_\eta |C(0)|,
\end{align}
which is finite by \eqref{e2.5}. Another consequence of
Proposition \eqref{prop2.1} is
\begin{align}\label{e2.127}
 \gamma_\xi=\gamma_\eta\, \E^0_\eta |C(0)|^{-1}.
\end{align}
The number $\E^0_\eta |C(0)|^{-1}$ might be called the number of
clusters per vertex in percolation theory; see e.g.\ \cite{Grimmett99}.

We also need the following consequence of Proposition \eqref{prop2.1}.

\begin{lemma}\label{lemma1} The relationship \eqref{nozeroeta} is equivalent to
\begin{align}\label{nozero}
\P^0_{\xi}\bigg(\sum_{x\in D(0)} e^{\mathbf{i}wx}\ne 0\bigg)>0,\quad w\in\R^d.
\end{align}
\end{lemma}
\emph{Proof:} The relationship \eqref{nozero} holds iff
\begin{align*}
\E^0_{\xi}\bigg|\sum_{x\in D(0)} e^{\mathbf{i}wx}\bigg| \neq 0.
\end{align*}
By \eqref{e2.15} this is equivalent with
\begin{align*}
0\neq\E^0_{\eta}\bigg|\sum_{x\in C(0)-\tau(0)} e^{\mathbf{i}wx}\bigg|
=\E^0_{\eta}\bigg|e^{-\mathbf{i}w\tau(0)}\sum_{x\in C(0)} e^{\mathbf{i}wx}\bigg|
=\E^0_{\eta}\bigg|\sum_{x\in C(0)} e^{\mathbf{i}wx}\bigg|.
\end{align*}
This implies the assertion.\qed

\bigskip

The Fourier transform of $P$ is the function $\hat{P}\colon \R^d\to\mathbb{C}$ given by
\begin{align*}
\hat{P}(w):=\int P(x)e^{\mathbf{i}wx}\,dx,\quad w\in\R^d.
\end{align*}
This transform can be expressed in terms of the typical cluster:

\begin{proposition}\label{prop2.3} For any $w\in\R^d$,
\begin{align}\label{e2.17}
\gamma_\eta+\gamma_\eta^2\hat{P}(w)
=\gamma_\xi\, \E^0_{\xi}\bigg|\sum_{x\in D(0)} e^{\mathbf{i}wx}\bigg|^2.
\end{align}
\end{proposition}
\emph{Proof:} First we apply \eqref{A5} with $(\omega, x) \mapsto
\ind\{x \in C(\omega, 0)\} e^{\mathbf{i} w y}$ to obtain that
\begin{align*}
\gamma_\eta+\gamma_\eta^2\hat{P}(w)
&=\gamma_\eta+\gamma_\eta^2\int \P^{0,x}(x\in C(0))e^{\mathbf{i}wx}\rho(x)\,dx\\
&=\gamma_\eta\,\E^0_\eta\sum_{x\in C(0)}e^{\mathbf{i}wx},
\end{align*}
where we recall the integrability assumption \eqref{e2.5}.

Since the clusters exhaust the points of $\eta$ we obtain that
\begin{align*}
\gamma_\eta+\gamma_\eta^2\hat{P}(w)
=\gamma_\eta\,\E^0_\eta\sum_{x\in\xi}\sum_{y\in C(0)}e^{\mathbf{i}wy}\I\{\tau(y)=x\}.
\end{align*}
Using the exchange formula \eqref{A2} (to be justified below) we get
\begin{align*}
\gamma_\eta+\gamma_\eta^2\hat{P}(w)
=&\gamma_\xi\,\E^0_{\xi}\sum_{x\in\eta} \sum_{y\in C(\theta_x,0)}e^{\mathbf{i}wy}
\I\{\tau(\theta_x,y)=-x\}\\
=&\gamma_\xi\,\E^0_{\xi}\sum_{x\in\eta}\sum_{y\in\eta} \I\{y-x\in C(\theta_x,0)\}e^{\mathbf{i}w(y-x)}
\I\{\tau(\theta_x,y-x)=-x\}\\
=&\gamma_\xi\,\E^0_{\xi}\sum_{x\in\eta} \sum_{y\in C(x)}e^{-\mathbf{i}wx}e^{\mathbf{i}wy}
\I\{\tau(y)=0\},
\end{align*}
where we have used the invariance properties
\eqref{covariance} and \eqref{e2.4}. For $0\in\xi$ and $x,y\in \eta$ the relations
$y\in C(x)$ and $\tau(y)=0$ are equivalent with $x,y\in D(0)$. Hence
\begin{align*}
\gamma_\eta+\gamma_\eta^2\hat{P}(w)=
\gamma_\xi\,\E^0_{\xi}\sum_{x\in D(0)}e^{-\mathbf{i}wx}\sum_{y\in D(0)}e^{\mathbf{i}wy},
\end{align*}
implying (by Fubini's theorem) the asserted formula \eqref{e2.17}.
The use of both the exchange formula and Fubini's theorem is justified by
$\E^0_{\xi}|D(0)|^2<\infty$, a consequence of \eqref{e2.12}
and assumption \eqref{e2.5}.\qed

\bigskip
\emph{Proof of Theorem \ref{tOZ}:} We shall use a classical theorem by Wiener
on the inversion of Fourier transforms. Recall that a {\em finite signed measure} $\mu$ on $\R^d$
is the difference of two finite measures. The Fourier transform of such a $\mu$
is defined by
\begin{align*}
\hat{\mu}(w):=\int e^{\mathbf{i}wx}\,\mu(dx),\quad w\in\R^d.
\end{align*}
The convolution $\mu\ast\nu$ of two finite signed measures $\mu$ and $\nu$ is the
finite signed measure defined by
\begin{align*}
\mu\ast\nu(B):=\iint\I\{x+y\in B\}\,\mu(dx)\,\nu(dy),\quad B\in\mathcal{B}^d.
\end{align*}
Note $\mu\ast\delta_0=\mu$, where $\delta_0(B):=\I\{0\in B\}$, $B\in\cB^d$.
Also note that $\widehat{\mu\ast\nu}=\hat{\mu}\hat{\nu}$.
Each $f\in L^1$ defines a finite signed measure $\mu_f:=\int\I\{x\in\cdot\}f(x)\,dx$.
(Later we will abuse notation and write $f$ instead of $\mu_f$.)
For $f,g\in L^1$ we have $\mu_f\ast\mu_g=\mu_{f\ast g}$.

Let $M^1$ denote the vector space of all finite signed measures of the form
$r\delta_0+\mu_f$, where $r\in\R$ and $f\in L^1$.
The Ornstein-Zernike equation \eqref{OZ} can be written as
\begin{align}\label{OZ5}
\mu_P=t\mu_Q\ast \nu,
\end{align}
where $t:=\gamma_\eta$ and $\nu:=t^{-1}\delta_0+\mu_P\in M^1$.
Proposition \ref{prop2.3}, assumption \eqref{nozeroeta} and Lemma \ref{lemma1}
imply that $\hat\nu(w)\ne 0$ for all $w\in\R^d$.
A theorem of Wiener (see Satz 13.2 in \cite{Joergens70}) says that $\nu$ can be inverted
within the convolution algebra $M^1$. This means that there is an $f\in L^1$ such that
$$
\nu\ast(t\delta_0+t^2\mu_f)=\delta_0.
$$
The function $Q:=P+tP\ast f$ is in $L^1$. Moreover,
\begin{align*}
t\mu_Q\ast\nu&=t\nu\ast(\mu_P+t\mu_{P\ast f})=\nu\ast(t\mu_P+t^2\mu_P\ast \mu_f)\\
&=\nu\ast\mu_P\ast(t\delta_0+t^2 \mu_f)=\mu_P\ast\delta_0=\mu_P,
\end{align*}
as required by \eqref{OZ5}.

To show that $Q$ is bounded, we apply \eqref{young-ineq-inf} to obtain
\begin{align}
\|Q\|_\infty = \|P\|_\infty+\gamma_\eta \|Q \ast P\|_{\infty }
\le \|P\|_\infty+\gamma_\eta \|Q \|_1 \|P\|_{\infty } < \infty .
\end{align}
\qed

\bigskip

The solution $Q$ of the OZ-equation \eqref{OZ}
has good integrability properties and can be used
to express the mean size of the cluster containing a typical point:

\begin{proposition}\label{prop2.2} Under the assumption of Theorem \ref{tOZ}
we have that
\begin{align}\label{5.105}
0\le \int Q(x)dx< \gamma_\eta^{-1}
\end{align}
and
\begin{align}\label{5.10}
\E^0_\eta|C(0)|=\bigg(1-\gamma_\eta\int Q(x)dx\bigg)^{-1}.
\end{align}
\end{proposition}
\emph{Proof:} Equation \eqref{e2-12} and the OZ-equation \eqref{OZ}
imply that
\begin{align*}
\E^0_\eta|C(0)|&=1+\gamma_\eta\int Q(x)dx+\gamma^2_\eta\int P(x)dx\int Q(x)dx\\
&=1+\gamma_\eta\int Q(x)dx+\gamma_\eta\big(\E^0_\eta|C(0)|-1\big)\int Q(x)dx.
\end{align*}
It follows that
\begin{align}\label{5.8}
\E^0_\eta|C(0)|=1+\gamma_\eta\E^0_\eta|C(0)|\int Q(x)dx.
\end{align}
Since $\E^0_\eta|C(0)|\ge 1$ we conclude that
$\int Q(x)dx\ge 0$. Moreover, since $\E^0_\eta|C(0)|<\infty$
we have $\gamma_\eta\int Q(x)dx<1$ and hence \eqref{5.10}.
\qed

\bigskip

\section{The random connection model}\label{secrcm}

In this section we consider
a stationary Poisson process $\eta_t$ on $\R^d$ with intensity $t\ge 0$
together with a measurable function $\varphi\colon\R^d\rightarrow[0,1]$  satisfying
\begin{align}\label{e19.1}
\varphi(x)=\varphi(-x),\quad x\in\R^d
\end{align}
and
\begin{align}\label{varphiint}
\int\varphi(x)\,dx<\infty.
\end{align}
Suppose any two distinct points $x,y\in\eta_t$ are adjacent with
probability $\varphi(y-x)$ independently of all other pairs and
independently of $\eta_t$. This yields the random connection model
(RCM), an undirected {\em random graph} $\cG$ with vertex set
$\eta_t$.  Each $x\in\eta_t$ belongs to a uniquely defined connected
component $C'(x)$.  The mapping $\tau$ from Section \ref{secOZ} is
defined as follows.  If $x\in\eta_t$ and $|C'(x)|<\infty$ then we let
$\tau(x)$ be the lexicographic minimum of $C'(x)$. (For all other
$x\in\R^d$ we let $\tau(x):=x$.) Hence, if all connected components of
$\cG$ are finite, the set of clusters consists exactly of these
connected components.

The Gilbert graph (briefly discussed in the introduction) based
on $\eta_t$ and a {\em gauge body} $B\subset\R^d$, that is a compact
and connected set containing the origin $0\in\R^d$, is a
special case of the RCM. It is obtained by choosing
$$
\varphi(x)=\I\{(B+x)\cap B\ne\emptyset\}.
$$
In contrast to the RCM, the Gilbert graph contains no additional randomness. Two points
$x,y\in\eta_t$ are adjacent if the shifted gauge bodies $B+x$ and
$B+y$ overlap.

In the next sections we shall study the properties of the
pair-connectedness function $P_t$ of the RCM and the solution $Q_t$ of the
associated OZE. In particular we shall show  that $P_t$
and $Q_t$ are analytic and relate the coefficients of their series
representation at $0$.  To do this properly we need to introduce the
model in a more formal way.  If the intensity $t$ is positive, then
$\eta_t$ can be (almost surely) represented as
\begin{align*}
  \eta_t=\{X_i:i\in\N\},
\end{align*}
where the $X_i$, $i\in\N$, are a.s.\ distinct random elements in $\R^d$.
For $t=0$ the Poisson process $\eta_t$ has (almost surely) no point.
Let $\R^{[2d]}$ denote the
space  of all sets $e\subset\R^d$ containing exactly two
elements. Any $e\in \R^{[2d]}$ is a potential edge of the RCM.
When equipped with the Hausdorff metric (see \cite{SW08}) this space is
a Borel subset of a complete separable metric space.
Let $<$ denote the strict lexicographic ordering on $\R^d$.
Introduce independent random variables $U_{i,j}$,
$i,j\in\N$, uniformly distributed on the unit interval $[0,1]$ such that
the double sequence $(U_{i,j})$ is independent of $\eta_t$.
For $t>0$
\begin{align}\label{e2.9}
\chi_t:=\{(\{X_i,X_j\},U_{i,j}):X_i<X_j,i,j\in\N\},\quad t>0,
\end{align}
is a point process on $\R^{[2d]}\times[0,1]$.
For $t=0$ we let $\chi_t$ equal the zero measure.
Note that $\eta_t$ can be recovered from $\chi_t$.
For $t>0$ we can define the RCM as a deterministic functional
of $\chi_t$ by taking for $i\ne j$
and $X_i<X_j$ the set $\{X_i,X_j\}$ as an edge of $\cG$
iff $U_{i,j}\le\varphi(X_i-X_j)$.

Justified by assumption \eqref{e19.1} we
can introduce a measurable function $\varphi^*\colon \R^{[2d]}\rightarrow[0,1]$ by
\begin{align*}
\varphi^*(e):=\varphi(y-x),\quad e=\{x,y\}\in \R^{[2d]}.
\end{align*}
If $\tilde{\chi}$ is a point process on $\R^{[2d]}\times[0,1]$, we can define
a graph $\cG(\tilde\chi):=\cG(\tilde\chi)=(V(\tilde\chi),E(\tilde\chi))$ as follows.
The vertex set is given by
\begin{align}\label{e2.37}
V(\tilde\chi):=\{e^-,e^+: e\in \R^{[2d]},\tilde{\chi}(\{e\}\times[0,1])=1\},
\end{align}
where $e^-$ and $e^+$ are the endpoints of $e\in\R^{[2d]}$.
A set $e\in \R^{[2d]}$ belongs to the edge set $E(\tilde\chi)$
of this graph iff $\tilde\chi(\{(e,u)\})=1$
for some $u\in[0,1]$ with $u\le\varphi^*(e)$.
In this notation our RCM is given as $\cG(\chi_t)$.
(For $t=0$ this is the empty graph.)
For $x\in V(\tilde\chi)$ we denote the cluster of $x$
(the connected component of $\cG(\tilde\chi)$) by $C(x,\tilde\chi)$.
(For convenience we set $C(x,\tilde\chi):=\{x\}$ for all other
$x\in\R^d$.)

In the remaining part of this section we give a few fundamental
results on the RCM that will be needed later but cannot be found in the literature.
We extend the (double) sequence $(U_{i,j})_{i,j=1}^\infty$ featuring in \eqref{e2.9}
to a sequence $(U_{i,j})_{i,j=0}^\infty$ of independent random variables uniformly
distributed on $[0,1]$, independent of the Poisson process $\eta_t$.
For $t>0$ we then define a
point process $\chi^0_t$ on $\R^{[2d]}\times[0,1]$ by
\begin{align}\label{chi0}
\chi^0_t:=\{(\{X_i,X_j\},U_{i,j}):X_i<X_j,i,j\in\N_0\},
\end{align}
where $\N_0:=\N\cup\{0\}$ and $X_0:=0$. The graph $\cG\big(\chi^0_t\big)$ can be interpreted as the
the RCM as seen from a {\em typical vertex} positioned at the origin.
For $x\in\R^d$ we define
\begin{align}\label{chi0x}
\chi^{0,x}_t:=\{(\{X_i,X_j\},U_{i,j}):X_i<X_j,i,j\in\N_{-1}\},
\end{align}
where $\N_{-1}:=\N_0\cup\{-1\}$, $X_{-1}:=x$ and $(U_{i,j})_{i,j=-1}^\infty$
has similar properties as $(U_{i,j})_{i,j=0}^\infty$.
In the case $t=0$ the point processes $\chi^0_t$
and $\chi^{0,x}_t$ are defined to be the empty set.

For $k \in \N$ we let $[k]:=\{1,2 ,\dots, k\}$. For any $x_1,\dots,x_k\in\R^d$
we introduce a random graph $\Gamma(x_1,\dots,x_k)$
with vertex set $\{x_1,\dots,x_k\}$ by taking
independent random  variables $U_{i,j}$, $i,j\in[k]$, with the
uniform distribution on
$[0,1]$  and by taking $\{x_i,x_j\}$  as an edge if
$x_i<x_j$ and $U_{i,j}\le \varphi(x_i-x_j)$. This is just the RCM
with a finite deterministic vertex set.
The next result is a version of Proposition 6.2 in \cite{MeesterRoy96}.
For the convenience of the reader we give a short proof.

\begin{proposition}\label{tclustersize} Let $n\in\N_0$ and set $x_0:=0$. Then
\begin{align}\label{3.41}\notag
\P\big(\big|C\big(0,\chi^0_t\big)\big|&=n+1\big)
=\frac{t^n}{n!}\int\P(\text{$\Gamma(x_0,\ldots,x_n)$ is connected})\\
&\quad \times\exp\bigg[-t \int \bigg(1-\prod^n_{i=0}(1-\varphi(y-x_i))\bigg)\,dy\bigg]
\,d(x_1,\dots,x_n).
\end{align}
In the case $n=0$ the right hand side has to be read as $\exp\big( -t \int\varphi(y)\,dy \big)$.
\end{proposition}
\emph{Proof:} We assume that $n\ge 1$. (The case $n=0$ is trivial.)
We have that $|C(0,\chi^0_t)|=n+1$ iff there are
$n$ distinct points $x_1,\dots,x_n\in\eta_t$ such that $\cG(\chi_t^0)$ restricted to those points
is connected and none
of the $x_i$ is connected to a point in $\eta_t\setminus\{x_1,\ldots,x_n\}$.
Given $\eta_t$, these two
events are (conditionally) independent and have
respective probabilities $\P(\Gamma(x_0,\dots,x_n) \text{ is connected})$ and
$$
\prod_{y\in\eta_t\setminus\{x_1,\ldots,x_n\}}\prod^n_{i=0}(1-\varphi(y-x_i)).
$$
After conditioning we obtain from the multivariate Mecke equation \eqref{Mecke}
that
\begin{align*}
&\P\big(\big|C\big(0,\chi^0_t\big)\big|=n+1\big)\\
&=\frac{t^n}{n!} \int \P(\text{$\Gamma(x_0,\dots,x_n)$ is connected})
\E \prod_{y\in\eta_t}\prod^n_{i=0}(1-\varphi(y-x_i)) \,d(x_1,\dots,x_n).
\end{align*}
Using the well-known formula for the characteristic functional
of $\eta_t$ (see e.g.\ \cite[Chapter 3]{LastPenrose16}) we get
the asserted formula \eqref{3.41}.
\qed

\bigskip

Next we need to discuss a Margulis-Russo type formula for
$\chi^{0,x}_t$. This formula provdides a power series expansion
of expectations of functions of $\chi^{0,x}_t$.
Adding just two points $0,x$ is enough for our purposes.
It would be no problem to extend the
result to a random connection model with any fixed number
of points added.
Let $n\in\N$ and
$\N_{-n-1}:=\N\cup\{0,-1,\ldots,-n-1\}$.  Extend the (double) sequence
$(U_{i,j})_{i,j=1}^\infty$ featuring in \eqref{e2.3} to a sequence
$(U_{i,j})_{i,j\in\N_{-n-1}}$ of independent random variables
uniformly distributed on $[0,1]$, independent of the Poisson process
$\eta_t$.  Let $x_0,x_{n+1}\in\R^d$ and $\bx=(x_1,\ldots,x_n)\in(\R^d)^n$
For $J\subset[n]$ we define $\bx_J:=(x_i)_{i\in J}$ and
\begin{align}\label{chi123}
\chi^{x_0,\bx_J,x_{n+1}}_t:=\{(\{X_i,X_j\},U_{i,j}):X_i<X_j,i,j\in \N_{J}\},
\end{align}
where $\N_J:=\N\cup\{0,-n-1\}\cup \{-i:i\in J\}$
and $(X_0,\ldots,X_{-n-1}):=(x_0,\ldots,x_{n+1})$.
In the case $t=0$
the indices $i,j$ are restricted to $\{-i:i\in J\}$.
Similarly we define the point process $\chi^{x_0,\bx_J}_t$.
For $J=\emptyset$ we set $\chi^{x_0}_t:=\chi^{x_0,\bx_\emptyset}_t$.
For $x_{n+1}:=x$ and $J=\emptyset$ the point process
$\chi^{0,\bx_{\emptyset },x}_t$ has the same distribution
as $\chi^{0,x}_t$ given by \eqref{chi0x}.
Let $f\colon\bN(\R^{[2d]}\times[0,1])\to\R$ be measurable and
fix some $x\in\R^d$. Define $F_t:=f\big(\chi^{0,x}_t \big)$ and
\begin{align}\label{Delta}
 \Delta^n_{\bx}F_t:=\sum_{J \subset [n]}(-1)^{n-|J|}
  f\big(\chi^{0,\bx_J,x}_t\big),\quad \bx\in(\R^d)^n.
\end{align}

We say that $f\colon\bN(\R^{[2d]}\times[0,1])\to\R$ is {\em determined} by
a compact subset $W\subset\R^d$, if $f(\mu)=f(\mu_W)$
for all $\mu\in\bN(\R^{[2d]}\times[0,1])$,
where
\begin{align}\label{3.67}
\mu_W:=\{(e,u)\in\mu:e\subset W\},
\end{align}
i.e.\ if the value of $f$ only depends on the edges with endpoints in
$W$.

\begin{theorem}\label{tperturbation} Let $f\colon \bN(\R^{[2d]}\times[0,1]) \to \R$
be measurable and let $x\in\R^d$.
Assume that $f$ is determined by a compact set $W\subset\R^d$
with $\{0,x\}\subset W$. Let
$s\ge 0$ and $t\ge -s$ such that $\E|F_{s+|t|}|<\infty$, where
$F_t:=f\big(\chi^{0,x}_t\big)$. Then
\begin{align}\label{Margulis2}
  \E F_{s+t}=\E F_s
  +\sum^\infty_{n=1}\frac{t^n}{n!}\int_{W^n}\E \Delta^n_{\bx}F_s\,d\bx.
\end{align}
\end{theorem}
\emph{Proof:} First we recall the Poisson process analogue of the Margulis-Russo
formula to be found in \cite{MolZu00}
and for a general phase space and more general integrability assumptions
in~\cite{Last14,LastPenrose16}.
Let $f\colon\bN(\R^d)\rightarrow\R$ be measurable, $n\in\N$ and $\bx=(x_1,\ldots,x_n)\in(\R^d)^n$.
Then we define a measurable function $D^n_{\bx}f\colon\bN(\R^d)\to\R$
by
\begin{align}\label{differenceop}
  D^n_{\bx}f(\mu):= \sum_{J \subset[n]}(-1)^{n-|J|}f\big(\mu\cup\{x_j:j\in J\}\big).
\end{align}
Assume now that there is a compact set $W\subset\R^d$ such that
$f(\mu)$ depends for each $\mu\in\bN(\R^d)$ only on the restriction
of $\mu$ to $W$. Then we have for all $s\ge 0$ and
$t\ge -s$ that
\begin{align}\label{Margulis}
  \E f(\eta_{s+t})=\E f(\eta_s)
  +\sum^\infty_{n=1}\frac{t^n}{n!}\int_{W^n} \big( \E D^n_{\bx}f(\eta_s) \big)\,d\bx,
\end{align}
provided that $\E|f(\eta_{s+|t|})|<\infty$.

For $\mu\in \bN(\R^d)$ and $x_1,\ldots,x_k\in\R^d$, $k\in\N$, we define
$\mu^{x_1,\ldots,x_k}:=\mu\cup\{x_1,\ldots,x_k\}$.
There is a probability kernel $K$ from $\bN(\R^d)$ to $\bN(\R^{[2d]}\times[0,1])$
such that for all $r\ge 0$
\begin{align}\label{3.12}
\P((\eta_r,\chi_r)\in\cdot)=\E\int \I\{(\eta_r,\psi)\in\cdot\}\,K(\eta_r,d\psi),
\end{align}
and, for any $x_0\in\R^d$, $n\ge 0$, and $\bx=(x_1,\ldots,x_n)\in(\R^d)^{n}$
\begin{align}\label{3.13}
\P\big(\big(\eta_t^{x_0,\bx},\chi^{x_0,\bx}_t \big)\in\cdot \big)
=\E\int \I\big\{\big(\eta_t^{x_0,\bx},\psi \big)\in\cdot \big\}
\,K\big(\eta_t^{x_0,\bx},d\psi \big).
\end{align}
Define a measurable function $f^*\colon \bN(\R^d)\to \R$ by
\begin{align*}
f^*(\mu):=\int f(\psi)\,K(\mu^{0,x},d\psi).
\end{align*}
By the triangle inequality and \eqref{3.13} for $n=0$,
\begin{align*}
\E |f^*(\eta_{s+|t|})|\le \E \int |f(\psi)|\,K(\eta^{0,x}_{s+|t|},d\psi)=\E |F_{s+|t|}|<\infty.
\end{align*}
The properties of the kernel $K$ imply that
\begin{align*}
f^*(\mu_W)&=\int f(\psi_W)\,K(\mu^{0,x}_W,d\psi)
=\int f(\psi_W)\,K(\mu^{0,x},d\psi)\\
&=\int f(\psi)\,K(\mu^{0,x},d\psi)=f^*(\mu).
\end{align*}
We can now apply \eqref{Margulis} with $f^*$ to obtain that
\begin{align}\label{Margulis9}
  \E f^*(\eta_{s+t})=\E f^*(\eta_s)
  +\sum^\infty_{n=1}\frac{t^n}{n!}\int_{W^n}\E D^n_{\bx}f^*(\eta_s)\,d\bx.
\end{align}
By \eqref{3.13} we have $\E f^*(\eta_{s+t})=\E F_{s+t}$ and $\E f^*(\eta_s)=\E F_s$.
Furthermore,
\begin{align*}
  \E D^n_\bx f^*(\eta_s)
&=\sum_{J \subset [n]}(-1)^{n-|J|}
 \E f^*(\eta_s\cup\{x_j:j\in J\})\\
&=\sum_{J \subset [n]}(-1)^{n-|J|}
 \E \int f(\psi)\, K\big(\eta^{0,x}_s\cup\{x_j:j\in J\},d\psi\big)\\
&=\sum_{J \subset [n]}(-1)^{n-|J|}
 \E f\big(\chi^{0,\bx_J,x}_s\big).
\end{align*}
In view of the definition \eqref{Delta}, we obtain the assertion.
\qed

\bigskip

We also need the following version of Proposition \ref{tclustersize}.
The proof is omitted.

\begin{proposition}\label{tclustersize2} Let $n\in\N_0$ and $x\in\R^d$. Then
\begin{align*}
\P\big(x\in C\big(0,\chi^{0,x}_t\big),& \big|C(0,\chi^{0,x}_t)\big|=n+2\big)
=\frac{t^n}{n!}\int\P(\text{$\Gamma(x_0,\ldots,x_{n+1})$ is connected})\\
&\quad \times\exp\bigg[-t \int \bigg(1-\prod^{n+1}_{i=0}(1-\varphi(y-x_i))\bigg)\,dy\bigg]
\,d(x_1,\dots,x_n),
\end{align*}
where $x_0:=0$ and $x_{n+1}:=x$.
\end{proposition}

\section{The OZE for the random connection model}\label{secOZErcm}

In this section we consider Poisson processes $\eta_t$ with intensity $t\ge 0$
and the associated RCM $\cG(\chi_t)$ as introduced
in the previous section. We assume that
\begin{align}\label{4.1}
0<m_\varphi<\infty,
\end{align}
where $m_\varphi:=\int \varphi(x)\,dx$.

The {\em critical intensity} is given by
$$
t_c:=\sup\big\{ t\ge 0: \P\big( \big| C\big( 0,\chi^0_t \big) \big|<\infty \big)=0 \big\}.
$$
For the Gilbert graph (in fact for general Boolean models)
it was proved in \cite{Hall88} that $0<t_c<\infty$.
The same is true for the more general RCM; see \cite[Theorem 6.1]{MeesterRoy96}.
By \cite[Theorem 6.2]{MeesterRoy96} we have that
\begin{align}\label{4.3}
t_c=\sup\big\{ t\ge 0: \E \big| C\big( 0,\chi^0_t \big) \big|<\infty \big\}.
\end{align}
We need to consider another critical intensity, namely
\begin{align}\label{4.4}
t_*:=\sup\big\{ t\ge 0:\text{$\E\exp\big(z\big| C\big( 0,\chi^0_t \big) \big|\big)<\infty$ for some $z>0$} \big\}.
\end{align}
Clearly we have $t_*\le t_c$.
For the Gilbert graph it is basically well-known that
$t_*= t_c$. (The arguments from \cite{Pen03} can be extended
from a convex and symmetric gauge body to a general $B$.)
We are not aware of a similar result for the RCM.
However, one can show, that
\begin{align}\label{4.5}
t_*\ge m^{-1}_\varphi.
\end{align}
This is due to the fact, that for $t < m^{-1}_\varphi$ the number of points in the cluster of the origin
can be dominated by the total progeny of a subcritical Galton-Watson process with
Poisson offspring distribution with mean $t m_{\varphi} < 1$; see the proof of
Theorem 6.1 in \cite{MeesterRoy96}.
It is well known, that this progeny is Borel distributed and hence has exponential moments \cite{Otter49}. 

By  \eqref{Mecke} the pair correlation function $\rho_t$ of $\eta_t$
satisfies $\rho_t\equiv 1$, so that the two-point Palm probability measures $\P^{x,y}_{\eta_t}$
of $\eta_t$ are well-defined. They are given  given by
the following lemma. Recall the definition \eqref{chi0} of $\chi^0_t$ and the
definition \eqref{chi123} of $\chi^{x,y}_t$.

\begin{lemma}\label{l4.1} We have
$\P^0_{\eta_t}\big(\chi_t\in\cdot\big)=\P\big(\chi^0_t\in\cdot\big)$.
Moreover, the Palm probability $\P^{x,y}_{\eta_t}$ can be chosen such that
\begin{align}\label{4.15}
\P^{x,y}_{\eta_t}\big(\chi_t\in\cdot\big)&=\P\big(\chi^{x,y}_t\in\cdot\big),\quad x,y\in\R^d.
\end{align}
\end{lemma}
\emph{Proof:} We prove the second formula. Let
$f\colon \bN(\R^{[2d]}\times[0,1])\times\R^d\times\R^d\to[0,\infty)$ be measurable.
Then we obtain from \eqref{3.12} and the Mecke equation \eqref{Mecke} that
\begin{align*}
\E\sumneq_{x,y\in\eta_t} f(\chi_t,x,y)
&=\E\sumneq_{x,y\in\eta_t}\int f(\psi,x,y)\,K(\eta_t,d\psi)\\
&=t^2\,\E \iint f(\psi,x,y)\,K\big(\eta^{x,y}_t,d\psi\big)\,d(x,y)\\
&=t^2\,\E \int f\big(\chi^{x,y}_t,x,y\big)\,d(x,y),
\end{align*}
where $\sum^{\ne}$ denotes summation over all pairs of distinct elements
of $\eta$ (a notation that is also used for multi-indices) and where
we have used \eqref{3.13} to get the final identity.
Comparing this with \eqref{A2Palm} (and using
that the pair correlation function $\rho_t$ of $\eta_t$
satisfies $\rho_t\equiv 1$), shows that
\eqref{4.15} holds for almost every $(x,y)$ (with respect to
Lebesgue measure on $\R^d\times\R^d$.
This shows the assertion.\qed

\bigskip

By Lemma \ref{l4.1} the pair connectedness function $P_t$ of the
RCM $\cG(\chi_t)$ is given by
$P_t(x,y)=P_t(y-x)$, where
\begin{align}
P_t(x):=\P\big(x\in C\big(0,\chi^{0,x}_t\big)\big).
\end{align}

\begin{theorem}\label{tOZRCM} Let $t<t_c$. Then there is a unique $Q_t\in L^1\cap L^\infty$ such that
\begin{align}\label{OZGilbert}
P_t=Q_t+t Q_t\ast P_t.
\end{align}
\end{theorem}
\emph{Proof:} We wish to apply Theorem \ref{tOZ}.
For any $x\in\R^d$ we define $\tau(x):=x$ if $x$ is not a member of
a finite cluster in $\cG(\chi_t)$. Otherwise we define $\tau(x)$ as the lexicographic
minimum of the cluster $C(x,\chi_t)$. Then
we have almost surely that $C(x)=C(x,\chi_t)$ for all $x\in\eta_t$,
where $C(x)$ is given by \eqref{e2.3}.
Since $t<t_c$ the integrability assumption \eqref{e2.5} follows from
\eqref{4.3}. Since the factorial moments measures of $\eta_t$ coincide
with Lebesgue measure (see \eqref{Mecke}),
Assumption \eqref{nozeroeta} follows from
Remark \ref{rassumption}.
%
\qed

\bigskip

\begin{proposition}\label{prop37} We have
\begin{align}\label{4.105}
0\le \int Q_t(x)dx< t^{-1},\quad 0<t<t_c,
\end{align}
and
\begin{align}\label{4.107}
 \E \big| C\big( 0,\chi^0_t \big) \big|=\bigg(1-t\int Q_t(x)dx\bigg)^{-1}.
\end{align}
\end{proposition}
\emph{Proof:} The two assertions follow from Proposition \ref{prop2.2}
and Lemma \ref{l4.1}.\qed

\bigskip

\begin{remark}\rm It is a fair conjecture that
$\lim_{t\uparrow t_c} \E\big| C\big( 0,\chi^0_t \big) \big|=\infty$,
but we have not found this in the literature. Under this hypothesis
\eqref{4.107} would show that
\begin{align}\label{4.109}
  t_c\lim_{t\to t_c-}\int Q_t(x)dx=1.
\end{align}
\end{remark}

In what follows we consider a measurable function
$g\colon\bN(\R^d)\to\R$ and fix some $x\in\R^d$. We study the function
$t \mapsto \E g\big(C\big(0,\chi^{0,x}_t\big)\big)$.
The results will imply, that
$t \mapsto P_t(x)$ and $t \mapsto Q_t(x)$ are analytic functions on $[0,t_*)$.
We assume that for any $\varepsilon>0$ there is an $n_0\in\N$, such that
\begin{align}\label{growthbound}
  |g(\mu)| \le e^{\varepsilon \mu(\R^d)}, \quad \mu \in \bN(\R^d),\ \mu(\R^d) \ge n_0.
\end{align}
This implies, that there is an $n_1\in\N$ such that
\begin{align}\label{growthbound2}
|g(\mu)| \le \exp\big( \varepsilon \max\big\{ \mu(\R^d), n_1 \big\} \big).
\end{align}

\begin{theorem}\label{analyticityOfClusterFctls} Suppose that $g\colon\bN(\R^d)\to\R$
satisfies \eqref{growthbound} and let $x\in\R^d$. Then the function
$t\mapsto \E g\big(C\big(0,\chi^{0,x}_t\big)\big)$ is analytic on $[0,t_*)$.
The expansion at $s\in[0,t_*)$ is given by
\eqref{Margulis2} with $F_t:=g\big(C\big(0,\chi^{0,x}_t\big)\big)$.
\end{theorem}

For the proof of Theorem \ref{analyticityOfClusterFctls}
we derive some preliminary results, that might be of independent interest.
Let $s,t\in[0,t^*)$ and define
\begin{align}
G(x,t):=\E g\big(C\big(0,\chi^{0,x}_t\big)\big).
\end{align}
We take a compact set $W\subset\R^d$ with $\{0,x\}\subset W$ and approximate the
function $G$ with
\begin{align}
  G_W(x,t):=\E g\big(C\big(0,\chi^{0,x}_{t,W}\big)\big),
\end{align}
where $\chi^{0,x}_{t,W}:=\big( \chi^{0,x}_t \big)_W$; see \eqref{3.67}. Note
that $\mu\mapsto g(C(0,\mu_W))$ is determined by $W$.
Let $z>0$ such that $\E\exp\big(2z\big|C\big(0,\chi^0_t\big)\big|\big)<\infty$.
Since
\begin{align}\label{4.23}
\big|C\big(0,\chi^{0,x}_t\big)\big|\le |C\big(0,\chi^0_t\big)\big|+|C\big(x,\chi^x_t\big)\big|
\end{align}
and $C\big(x,\chi^x_t\big)\overset{d}{=}C\big(0,\chi^0_t\big)$ we have that
\begin{align*}
\E\exp\big(z\big|C\big(0,\chi^{0,x}_t\big)\big|\big)<\infty.
\end{align*}
Choosing $\varepsilon=z$ in assumption \eqref{growthbound2} we obtain that
\begin{align*}
\E \big|g\big(C\big(0,\chi^{0,x}_{t,W}\big)\big)\big|
&\le \E\exp\big(z\max\big\{\big|C\big(0,\chi^{0,x}_{t,W}\big)\big|,n_1\big\}\big) \\
&\le \E\exp\big(z\max\big\{\big|C\big(0,\chi^{0,x}_t\big)\big|,n_1\big\}\big)<\infty.
\end{align*}
Therefore we can apply Theorem \ref{tperturbation} to obtain that
\begin{align}\label{G_W(t)_series_representation}
  G_W(x,t) =\E g\big(C\big(0,\chi^{0,x}_{s,W}\big)\big)
+\sum_{n=1}^{\infty }(t-s)^n g_{W,n}(x,s),\quad t<t_*,
\end{align}
where
\begin{align}\label{def_g_n(t_0,W)}
g_{W,n}(x,s):=\frac{1}{n!}\int_{W^n} \E\bigg[\sum_{J \subset[n]} (-1)^{n-|J|}
g\Big( C\Big(0,\chi^{0,\bx_J,x}_{s,W}\Big)\Big)\bigg]\ d\bx,
\end{align}
and $\chi^{0,\bx_J,x}_{s,W}:=\big(\chi^{0,\bx_J,x}_{s}\big)_W$.
We use this definition for all Borel sets $W\subset\R^d$.

To bound the coefficients $g_{W,n}(x,t)$, we use the following
integral inequality. 
Recall
that $\Gamma(0,x_1,\ldots,x_n)$ denotes a RCM with vertex set
$\{0,x_1,\ldots,x_n\}$.

\begin{lemma}\label{le:integral-inequality}
Let $n\in\N$. Then
  \begin{align}
\int \P(\text{$\Gamma(0,x_1,\ldots,x_n)$ is connected})\,d(x_1,\dots,x_n)\le n! m_\varphi^ne^{n+1}.
  \end{align}
Furthermore we have for any $x\in\R^d$ that
  \begin{align}
\int \P(\text{$\Gamma(0,x_1,\ldots,x_n,x)$ is connected})\,d(x_1,\dots,x_n)\le n! m_\varphi^ne^{n+2}.
  \end{align}
\end{lemma}
\emph{Proof:} We prove the second inequality. For all $a_0,\ldots,a_{n+1}\in[0,1]$ we have the inequality
\begin{align*}
1-\prod^{n+1}_{i=0}(1-a_i)\le \sum^{n+1}_{i=0}a_i.
\end{align*}
For any $t>0$ we therefore obtain from Proposition \ref{tclustersize2}  that
\begin{align*}
1&\ge \frac{t^n}{n!}\int\P(\text{$\Gamma(x_0,\ldots,x_n,x)$ is connected})
\exp\bigg( -t \sum^{n+1}_{i=0}\int \varphi(y-x_i)\,dy \bigg)\,d(x_1,\dots,x_n)\\
&\ge \frac{t^n}{n!}\int\P(\text{$\Gamma(x_0,\ldots,x_n,x)$ is connected})
\exp(-t(n+2)m_\varphi)\,d(x_1,\dots,x_n).
\end{align*}
Choosing $t=m_{\varphi}^{-1}$ yields the asserted inequality.\qed

\bigskip

There is a qualitative difference between the study of analyticity of $G$ at $s=0$ and $s>0$.
In fact the condition \eqref{growthbound}
can be slightly relaxed for $s= 0$.

\begin{lemma}\label{le:bound_g_n(0,W)}
  Let $n \in \N$ and assume that there is a constant $c \geq 1$, such that
  \begin{align}\label{growthbound3}
    |g(\mu)| \leq c^{\mu(\R^d)}, \quad \mu \in \bN(\R^d),\ \mu(\R^d) \leq n+1.
  \end{align}
Let $x\in\R^d$ and $W\subset\R^d$ be a Borel set such that $\{0,x\}\subset W$. Then
  \begin{align}
    |g_{W,n}(0,x)|\le c(c+1)(1+e)\big((1+c)m_{\varphi}e\big)^n.
  \end{align}
\end{lemma}
\emph{Proof:} Let $x_0:=0$ and $x_{n+1}:=x$.
Take $\bx=(x_1,\ldots,x_n)\in W^n$.
We recall, that for $s=0$, the point process $\eta_s$ is the zero measure
and $\Gamma(x_0,\ldots,x_{n+1})=\cG\big(\chi^{x_0,\bx,x_{n+1}}_0\big)$
is the RCM with vertex set $\{x_0,\ldots,x_{n+1}\}$.
Let $i\in [n]$ such that $x_0$ and $x_i$ are not connected by a path in $\Gamma(x_0,\ldots,x_{n+1})$.
Then we have for any $J\subset[n]\setminus\{i\}$ that
\begin{align*}
  (-1)^{n-|J|}
  g\big(C\big(0,\chi^{x_0,\bx_J,x_{n+1}}_0\big)\big)+(-1)^{n-|J \cup \{i\}|}
  g\big(C\big(0,\chi^{x_0,\bx_{J\cup\{i\}}, x_{n+1}}_0\big)\big) = 0,
\end{align*}
since the cluster of $0$ is the same in both summands.
If, on the other hand,
$x_0$ and $x_i$ are  connected in $\Gamma(x_0,\ldots,x_{n+1})$
for each $i\in[n]$ then either $\Gamma(x_0,\ldots,x_n)$ is connected
(and $x$ is not connected to any of the points $x_0,\ldots,x_n$) or
$\Gamma(x_0,\ldots,x_{n+1})$ is connected.
Hence we have that
\begin{align*}
  |g_{W,n}(0,x)| \leq \frac{1}{n!}\E \int h(x_1,\ldots,x_n,x)
\sum_{J\subset[n]} \big|g\big(C\big(0,\chi^{x_0,\bx_J,x_n}_0\big)\big) \big|\,d\bx.
\end{align*}
where
$$
h(x_1,\ldots,x_n,x):=\I\{\text{$\Gamma(x_0,\ldots,x_n)$ is connected}\}
+\I\{\text{$\Gamma(x_0,\ldots,x_{n+1})$ is connected}\}.
$$
Our assumption \eqref{growthbound3} and the binomial formula imply that
\begin{align*}
  |g_{W,n}(0,x)| &\le
\frac{c(c+1)^n}{n!}\int h(x_1,\ldots,x_n,x) \, d(x_1,\dots, x_n).
\end{align*}
Using here Lemma \ref{le:integral-inequality} concludes the proof.
\qed

\bigskip

In the following it is convenient to introduce a function $c_V\colon [0, \infty)\to[0,\infty)$
with $c_V(t) > 0$ for $t < t_*$ satisfying
\begin{align}\label{exponential_decay_clustersize}
  \P\big( \big| C\big( 0,\chi^{0,x}_t \big) \big|=n \big) \le e^{-c_V(t) n},\quad n\in\N,\,x\in\R^d,\,\,t<t_*.
\end{align}
By definition of $t_*$ and \eqref{4.23}, such a function exists.
We next bound $g_n(t,W)$ for $t>0$.

\begin{lemma}\label{le:bound_g_n(t,W)}
Let $t \in (0,t_*)$ and $x\in\R^d$. Assume that \eqref{growthbound} holds. Then there is an $n_0(t)\in\N$
such that for all $n>n_0(t)$ and all Borel sets $W\subset\R^d$ with $\{0,x\}\subset W$
  \begin{align}
    |g_{W,n}(x,t)|\le \frac{e^{-c_V(t)/2}}{1-e^{-c_V(t)/2}}\left(\frac{2e^{-c_V(t)/2}}{t(1-e^{-c_V(t)/2})}\right)^n.
  \end{align}
\end{lemma}
\emph{Proof:} As before we set $x_0:=0$ and $x_{n+1}:=x$. Let $n\in\N$.
With the same argument as in the proof of Lemma \ref{le:bound_g_n(0,W)} we conclude,
that 
\begin{align}\label{eq:bound_g_n(t,W):1}\notag
  |g_{W,n}(x,t)| & \leq \frac{1}{n!} \int_{W^n}
  \E\bigg[\I\big\{\text{$x_i\in C\big(0,\chi^{x_0,\bx,x_{n+1}}_{t,W}\big)$
for all $\in [n]$}\big\}  \\
&\hspace{3cm} \times \sum_{J \subset [n]}
\big|g\big(C\big(0, \chi^{x_0,\bx_J,x_{n+1}}_{t,W} \big)\big)\big|\bigg] \,d\bx.
\end{align}
Setting $\varepsilon:=c_V(t)/2$ and using \eqref{growthbound}
we find an $n_0\in\N$ such that
\begin{align*}
\sum_{J \subset [n]}
\big|g\big(C\big(0,\chi^{x_0,\bx_J,x_{n+1}}_{t,W} \big)\big)\big|
  & \leq \sum_{J \subset [n]} \exp\Big( \frac{c_V(t)}{2}
\max\big\{n_0,\big|C\big(0, \chi^{x_0,\bx_J,x_{n+1}}_{t,W} \big)\big|\big\}\Big) \\
  & \le 2^n \exp\Big( \frac{c_V(t)}{2}
\max\big\{n_0,\big|C\big(0, \chi^{x_0,\bx,x_{n+1}}_{t,W} \big)\big|\big\}\Big).
\end{align*}
Inserting this in \eqref{eq:bound_g_n(t,W):1} and using \eqref{3.13} yields
\begin{align*}
  |g_{W,n}(x,t)| \le &\frac{2^n}{n!} \E\iint
\I\big\{\text{$x_i\in C\big(0,\psi_W\big)$
for all $i\in [n]$}\big\}\\
&\times\exp\Big(\frac{c_V(t)}{2}\max\{n_0,|C(0, \psi_W )|\}\Big)
\,K\big(\eta_t^{x_0,\bx,x_{n+1}},d\psi \big)\,d\bx.
\end{align*}
The Mecke equation \eqref{Mecke} gives
\begin{align*}
  |g_{W,n}(x,t)| &\le \frac{2^n}{t^nn!}
\E\sumneq_{x_1,\ldots,x_n\in\eta_t}
\int\I\big\{\text{$x_i\in C\big(0,\psi_W\big)$ for all $i\in [n]$}\big\}\\
&\quad \times\exp\Big(\frac{c_V(t)}{2}\max\big\{n_0,\big|C\big(0, \psi_W \big)\big|\big\}\Big)
\,K(\eta^{0,x}_t,d\psi)\\
&=\frac{2^n}{t^nn!} \E\sumneq_{x_1,\ldots,x_n\in\eta_t}
\I\big\{\text{$x_i\in C\big(0,\chi^{0,x}_{t,W}\big)$ for all $i \in [n]$}\big\}\\
&\qquad \times\exp\Big(\frac{c_V(t)}{2}
\max\big\{n_0,\big|C\big(0, \chi^{0,x}_{t,W} \big)\big|\big\}\Big),
\end{align*}
where we have used \eqref{3.13} to achieve the final identity.
Therefore
\begin{align*}
  |g_{W,n}(x,t)| &\le \frac{2^n}{t^nn!}\sum^\infty_{k=1} \E\I\big\{
  \big|C \big(0,\chi^{0,x}_t\big) \big|=k
  \big\}\exp\Big(\frac{c_V(t)}{2}
  \max\big\{n_0,\big|C\big(0, \chi^{0,x}_{t,W} \big)\big|\big\}\Big)\\
&\qquad\times \sumneq_{x_1,\ldots,x_n\in\eta_t}\I\big\{\text{$x_i\in C\big(0,\chi^{0,x}_t\big)$
    for all $i \in [n]$}\big\}.
\end{align*}
On the event $\big\{ \big| C\big(0,\chi^{0,x}_t) \big|=k \big\}$ the above integral simplifies to
$(k-1)\cdots (k-n)$. Hence
\begin{align*}
  |g_{W,n}(x,t)| & \leq\frac{2^n}{n! t^n} \sum_{k=n+1}^\infty e^{
    \frac{c_V(t)}{2} \max\{n_0, k\}} (k-1)\cdots(k-n)\P\big( \big|
  C\big( 0,\chi^{0,x}_t \big) \big|=k \big).
\end{align*}
Finally we apply \eqref{exponential_decay_clustersize} and use
the well known formula for the factorial moment of the geometric
distribution, to get for $n > n_0$ that
\begin{align*}
  |g_{W,n}(x,t)| & \leq \frac{2^n}{n! t^n} \sum_{k = n+1}^\infty
e^{-\frac{c_V(t)}{2} k} (k-1) \cdot \ldots \cdot (k-n) =
  \frac{e^{-c_V(t)/2}}{1 - e^{-c_V(t)/2}}\left( \frac{2
      e^{-c_V(t)/2}}{t(1-e^{-c_V(t)/2})} \right)^n,
\end{align*}
as asserted. \qed

\bigskip

\emph{Proof of Theorem \ref{analyticityOfClusterFctls}:}
Let $W_k$, $k\in\N$, be a sequence of compact sets with union $\R^d$.
Since $\P\big(\big|C(0,\chi^{0,x}_t\big)\big|<\infty)=1$ we have
\begin{align*}
\lim_{k\to\infty}g\big(C\big(0,\chi^{0,x}_{t,W_k}\big)\big)=g\big(C\big(0,\chi^{0,x}_t\big)\big),
\quad \P\text{-a.s.}
\end{align*}
By \eqref{growthbound2},
\begin{align*}
  \big|g\big(C\big(0,\chi^{0,x}_{t,W_k}\big)\big)\big|
&\le \exp\big( c_V(t)\max\big\{\big|C\big(0,\chi^{0,x}_{t,W_k}\big)\big|,n_1\}\big)\\
&\le \exp\big( c_V(t)\max\big\{\big|C\big(0,\chi^{0,x}_t,\big)\big|,n_1\}\big).
\end{align*}
It follows from \eqref{exponential_decay_clustersize}, that
\begin{align*}
\E\exp\big(c_V(t)\big|C\big(0, \chi^{0,x}_t\big)\big|\big)<\infty.
\end{align*}
Dominated convergence implies
\begin{align*}
\lim_{k\to\infty} G_{W_k}(x,t)=\lim_{k\to\infty}\E g\big(C\big(0,\chi^{0,x}_{t,W_k}\big)\big)
=\E g\big(C\big(0,\chi^{0,x}_t\big)\big)=G(x,t).
\end{align*}
Similarly, dominated convergence implies for any $n\in\N$ that
\begin{align*}
\lim_{k\to\infty} g_{W_k,n}(x,s)=g_{\R^d,n}(x,s).
\end{align*}
Now we use the series representation \eqref{G_W(t)_series_representation}
for $W=W_k$.
The Lemmas \ref{le:bound_g_n(0,W)} and \ref{le:bound_g_n(t,W)} allow
us to apply dominated convergence to obtain that
\begin{align*}
  \lim_{k\to\infty} G_{W_k}(t)=G(s)+\sum_{n=1}^\infty (t-s)^n \lim_{k\to\infty} g_{W_k,n}(x,s)
=G(s)+\sum_{n=1}^\infty (t-s)^ng_{\R^d,n}(x,s)
\end{align*}
holds for all $t$ in some open neighborhood of
$s\in[0,t_c)$. This completes the proof.\qed

\bigskip
We want to point out, that due to the relaxed growth bound of Lemma \ref{le:bound_g_n(0,W)}
in comparison to \eqref{growthbound}, any functional that grows exponentially
in the size of the cluster of the origin is analytic at least in $s=0$.
The Lemmas \ref{le:bound_g_n(0,W)}
and \ref{le:bound_g_n(t,W)} also give a lower bound for the radius of convergence
of the series representation of $G(t)$ which is rather small though.

Theorem \ref{analyticityOfClusterFctls} shows that the
pair connectedness function and the expected cluster size are analytic
functions on the whole interval $[0,t_*)$. In particular, given $x\in\R^d$,
every $s\in [0,t_*)$ has a neighbourhood $U(s)$ such that
\begin{align*}
  P_t(x)=\sum_{n=0}^\infty (t-s)^n p_n(x,s), \quad t \in U(s),
\end{align*}
where $p_0(x,s) := P_s(x)$ and, for $n\in\N$,
\begin{align}\label{eq:def_p_n}
  p_n(x,s):=\frac{1}{n!}
\int \E \sum_{J \subset [n]} (-1)^{n-|J|}
\I\big\{ x\in C\big( 0,\chi^{0,\bx_J,x}_s \big) \big\} \ d\bx.
\end{align}
We summarize the integrability properties of the coefficients $p_n$ in the following corollary.

\begin{corollary}\label{analyticity-P_t}
  For any $n \in \N_0$ and $t \in [0, t_*)$ there are constants $c_1(t),c_2(t)$ such that
  \begin{align}
    \|p_n(\cdot, t)\|_{\infty } \leq c_1(t) c_2(t)^n, \\
    \|p_n(\cdot, t)\|_{1} \leq c_1(t) c_2(t)^n.
  \end{align}
  Moreover, for any $s \in [0, t_*)$ there is a neighbourhood $U(s)$ such that
  \begin{align}\label{eq:expansion_for_P_t}
    P_t( \cdot) = \sum_{n = 0}^\infty (t-s)^n p_n(\cdot, t), \quad t \in U(s),
  \end{align}
  where the convergence holds in $L^1$ and $L^\infty $.
\end{corollary}
\emph{Proof:} We observed in the proofs of Lemma
\ref{le:bound_g_n(0,W)} and \ref{le:bound_g_n(t,W)}, that the bounds
on $g_n$ only depend on the growth bounds \eqref{growthbound} or
\eqref{growthbound2} respectively, which immediately yields the bound
of $\|p_n(\cdot, t)\|_{\infty }$.

The arguments in the proof of Lemma \ref{le:bound_g_n(0,W)} show
\begin{align*}
  \int_{\R^d} |p_n(x, 0)| \ dx
\leq \frac{2^n}{n!} \int_{(\R^d)^{n+1}} \I\{\text{$\Gamma(0,x_1, \dots,x_n,x\}$ is connected}\}
\,d(x_1,\ldots,x_n,x)
\end{align*}
which can be bounded using Lemma \ref{le:integral-inequality}. The
bound on $\|p_n(\cdot,t)\|_{1}$ with $t > 0$ can be derived in a
similar way. It is clear, that these bounds imply the $L^1$ and
$L^\infty $ convergence of the sum in
\eqref{eq:expansion_for_P_t} for $t$ in a neighbourhood $U(s)$ of $s$.\qed

\bigskip

With a good understanding of the analyticity of $P_t$ we are now able
to show similar results for the solution $Q_t$ of the Ornstein-Zernike
equation. We will write $f^{\ast n}$ for an $n$-fold convolution of
the function $f$ with itself, i.e. $f^{\ast (n+1)} := f^{\ast n} \ast
f$ for $n \in \N$ and $f^{\ast 1} := f$. In the same spirit, we define
$\ast_{k = a}^b f_k:=f_b \ast (\ast_{k = a}^{b-1} f_k)$ and
$\ast_{k=a}^a f_k:=f_a$ for $a < b \in \Z$ and functions
$f_a,\ldots,f_b$.

\begin{proposition}\label{analyticity-Q_t}
  If $t \geq 0$ is, such that $\E\big| C\big( 0,\chi^0_t \big) \big|<2$, then
  \begin{align}\label{eq:anal-Q_t:1}
    Q_t = \sum_{n = 0}^\infty (-t)^n P_t^{\ast (n+1)}
  \end{align}
  in $L^1$ and $L^\infty $. Moreover, for any $s \in [0, t_*)$ there
  is a neighbourhood $U(s)$ and functions $q_n(\cdot,t_*)$, such
  that
  \begin{align}\label{eq:anal-Q_t:2}
    Q_t(\cdot)=\sum_{n=0}^\infty (t-s)^n q_n(\cdot,s), \quad t\in U(s)
  \end{align}
  in $L^1$ and $L^\infty $. The coefficients can
be recursively determined by the solvable equations
  \begin{align}
    q_{0}(\cdot, s)+s p_0(\cdot,s) \ast q_0(\cdot, s)
&= p_{0}(\cdot, s)\label{eq:analyticity_of_Q_t:rec1}, \\
    q_{n}(\cdot, s)+s q_{n}(\cdot, s)\ast p_{0}(\cdot, s)
&= p_{n}(\cdot, s) - \sum_{k=1}^{n} q_{n-k}(\cdot, s)\ast (p_{k-1}(\cdot, s)+s p_{k}(\cdot, s))
\label{eq:analyticity_of_Q_t:rec2}.
  \end{align}
\end{proposition}
\emph{Proof:} From $\E\big| C\big( 0,\chi^0_t \big) \big|<2$ and
\eqref{e2-12} we obtain $t \|P_t\|_1 < 1$. By
\eqref{young-ineq-1} and \eqref{young-ineq-inf} we have
\begin{align}
  \big\| P_t^{\ast k+1} \big\|_\infty  \leq \|P_t\|_1^k \|P_t\|_\infty, \quad k \in \N,
\end{align}
as well as
\begin{align}
  \big\| P_t^{\ast k+1} \big\|_1  \leq \|P_t\|_1^{k+1}, \quad k \in \N
\end{align}
and hence, the convergence of the right-hand-side of
\eqref{eq:anal-Q_t:1} in $L^1$ and $L^\infty $. A simple calculation
shows, that \eqref{eq:anal-Q_t:1} solves the Ornstein-Zernike
equation.

To prove the second part of the claim we start by solving
\eqref{eq:analyticity_of_Q_t:rec1} and
\eqref{eq:analyticity_of_Q_t:rec2} for $q_0(\cdot,s)$ and
$q_n(\cdot,s)$ respectively. From the proof of Theorem \ref{tOZ} we
know, that there is a function $f \in L^1$ such that $(\delta_0+s
p_0(\cdot,s)) \ast (\delta_0+s f) = \delta_0$. Hence the
equations \eqref{eq:analyticity_of_Q_t:rec1} and
\eqref{eq:analyticity_of_Q_t:rec2} are equivalent to
\begin{align}
  q_{0}(\cdot,s) &= (\delta_0+s f) \ast p_{0}(\cdot,s), \label{eq:anal-Q_t:3} \\
  q_{n}(\cdot,s) &= (\delta_0+s f) \ast \bigg(  p_{n}(\cdot,
    s) - \sum_{k=1}^{n} q_{n-k}(\cdot,s)\ast (p_{k-1}(\cdot,
    s)+s p_{k}(\cdot,s))  \bigg). \label{eq:anal-Q_t:4}
\end{align}
This implies that the $q_n$ can be recursively determined.

In the next step we show the series in
\eqref{eq:anal-Q_t:2} converges. We fix $s$ and write $p_n$ for
$p_n(\cdot,s)$ and $q_n$ for $q_n(\cdot,s)$. We choose $p,c\in\R$ such
that $\max\{\|p_n\|_1, \|p_n\|_\infty\} \leq p^n$ for all $n\in\N$
as well as $\max\{\|p_0\|_1, \|q_0\|_1, \|\delta_0+s f\|_1\}\leq c$.
This is possible due to Corollary \ref{analyticity-P_t} and
Theorem \ref{tOZ}. Moreover, we choose $q$ such that $q > p$,
\begin{align}
  q > c(p+c^2+s c p)
\end{align}
and
\begin{align}
  c \left( \frac{p}{q}+\frac{2c}{q}+\frac{p}{q(q-p)}+\frac{s c
      p}{q}+\frac{s p}{q - p}\right) \leq 1.
\end{align}
By \eqref{eq:anal-Q_t:4} and \eqref{young-ineq-1} we have
\begin{align}
\|q_1\|_1
=\left\|(\delta_0+s f)\ast(p_1-q_0\ast p_0 -sq_0\ast p_1) \right\|_1
\leq c(p+c^2+s c p)<q.
\end{align}
By induction over $n$ we obtain
\begin{align*}
  \|q_{n+1}\|_1
  & = \left\| (\delta_0+s f) \ast \left( p_{n+1} -
\sum_{k = 0}^n p_{k} \ast q_{n-k} - s \sum_{k=1}^{n+1} p_{k} \ast q_{n+1-k}\right) \right\|_1 \\
  & \leq c \left( p^{n+1}+p^n c+c q^n+\frac{p^n q - p q^n}{p-q}+s p^{n+1} c
+ s \frac{p^{n+1} q - p q^{n+1}}{p-q}\right) \\
  & = q^{n+1} c \left( \left( \frac{p}{q}\right)^{n+1}
+\left(\frac{p}{q}\right)^n
\frac{c}{q} +\frac{c}{q}+\frac{1}{q-p} \left( \left( \frac{p}{q}\right)
-\left( \frac{p}{q}\right)^n \right)+\right.\\
  &\hspace{2cm} \left.+s c \left( \frac{p}{q}\right)^{n+1}+s \frac{p}{q-p}
\left(1-\left( \frac{p}{q}\right)^n \right)\right).\\
  & \leq q^{n+1} c \left( \frac{p}{q}+\frac{2c}{q}+\frac{p}{q(q-p)}+
\frac{s c p}{q}+\frac{s p}{q - p}\right) \\
  & \leq q^{n+1}.
\end{align*}
If we use \eqref{young-ineq-inf} instead of \eqref{young-ineq-1} we
obtain the same bound for $\|q_n\|_\infty $ which implies the
convergence of the sum in \eqref{eq:anal-Q_t:2}.

It remains to show, that the sum in \eqref{eq:anal-Q_t:2} solves the
Ornstein-Zernike equation. This is achieved by rewriting the
Ornstein-Zernike equation in the form
\begin{align}
  P_t = Q_t+(t-s) P_t \ast Q_t+s P_t \ast Q_t.
\end{align}
Substituting for $P_t$ and $Q_t$ the series expansion at $s$ yields
that the equation holds, if for all $n \in \N_0$
\begin{align}
  p_n = q_n+\sum_{k=1}^n p_{k-1} \ast q_{n-k} +s\sum_{k=0}^np_k \ast q_{n-k}
\end{align}
which is equivalent to \eqref{eq:analyticity_of_Q_t:rec1} and \eqref{eq:analyticity_of_Q_t:rec2}.\qed

\bigskip

\section{Combinatorics for small intensities}\label{seccomb}

The coefficients $p_n$ in \eqref{eq:expansion_for_P_t} (given by \eqref{eq:def_p_n})
are quite complex probabilistic objects. In the
expansion of $P_t(x)$ around $s=0$ however, $\eta_{s}$ vanishes
and the only random objects that remain are the random connections
between the points $0,x_1,\dots,x_n,x$. This leads to an almost
combinatorial interpretation of the $p_n(x,0)$
(In the Boolean model with a fixed gauge body all randomness
disappears.)
Moreover this interpretation provides a simple combinatorial
way to determine the coefficients $q_n(x,0)$ in \eqref{eq:anal-Q_t:2}.

For $n \in \N_0$ let $\G_n$ be the set of connected graphs with $n+2$
vertices $\{0,\ldots,n+1\}$. For a graph
$G = (V(G),E(G))\in\G_n$ with vertex set $V(G)$ and edge set $E(G)$ we call $0$ the
start-vertex and $n+1$ the end-vertex.
For $i,j \in V(G)$ and $I \subset [n]$ we write "$i \leftrightarrow j$ in
$G|I$" if there is a path
from $i$ to $j$ in $G$ that uses only vertices in $I\cup\{i,j\}$.
For $n \in \N_0$ we define the combinatorial functionals
$\pi_n\colon\G_n \to \Z$ by
\begin{align}
  \pi_n(G):=\sum_{I \subset [n]} (-1)^{n - |I|}\ind\{0 \leftrightarrow n+1 \text{ in } G|I\}
\end{align}
By a slight abuse of notation we write $G = \cG\big( \chi_0^{x_0,
  \dots, x_{n+1}} \big)$ for $G \in \G_n$ if the two graphs are equal
after changing the labels in $G$ from $i$ to $x_i$.

It was shown in Lemma \ref{le:bound_g_n(0,W)} that the integrand in
\eqref{eq:def_p_n} vanishes if $\cG\big( \chi_0^{x_0,\bx,x_{n+1}}\big)$
is not connected. Hence \eqref{eq:def_p_n} is
equivalent to
\begin{align}
\begin{aligned}\label{eq:repr_p_n}
  & p_n(x,s)\\
  &=\frac{1}{n!} \int \E \sum_{G \in \G_n} \sum_{J \subset [n]}
  (-1)^{n-|J|} \ind\big\{ x \in C\big( 0, \chi_0^{0, \bx, x} \big) \big\}
  \ind\big\{ \cG\big( \chi^{0,\bx,x}_0 \big) = G\big\}\, d\bx\\
  &= \frac{1}{n!} \int \sum_{G \in \G_n} \pi_n(G)
\E\ind\big\{ \cG\big( \chi^{0,\bx,x}_0 \big) = G \big\}\, d\bx\\
  & = \frac{1}{n!} \sum_{G \in \G_n} \pi_n(G) I_n(G,x),
\end{aligned}
\end{align}
where $I_n:\G_n \times \R^d \to [0, \infty )$ is defined by
\begin{align}
\begin{aligned}
  I_n&(G,x):=\int \P\big(\cG\big( \chi^{0,\bx,x}_0\big)=G\big)\, d\bx\\
  &=\int \E \prod_{\{i,j\} \in E(G)}\ind\big\{ \{x_i,x_j\}\in E\big(\chi^{x_0,\bx,x_{n+1}}_0\big)\big\}
\prod_{\{i,j\} \notin E(G)} \ind\big\{ \{x_i,x_j\}
  \notin E\big( \chi^{x_0,\bx,x_{n+1}}_0\big)\big\}\,d\bx\\
  &=\int\prod_{\{i,j\} \in E(G)} \varphi(x_i-x_j)
  \prod_{\{i,j\} \notin E(G)} (1-\varphi(x_i-x_j)) \
  d(x_1,\ldots,x_n),
\end{aligned}
\end{align}
where again $x_0:=0$ and $x_{n+1}:=x$.  By \eqref{eq:repr_p_n} we have
found a representation of $p_n(x)$ as a sum over the graphs in $\G_n$
where each summand consists of a purely combinatorial factor and an
integral-geometric factor. This representation looks rather natural,
but is not well suited for the convolution. Therefore we will derive
a second representation, that convolutes in a very simple way. This
will also enable us, to give a very simple representation of the
$q_n(x)$.

Let $J_n\colon\G_n \times \R^d \to [0,\infty)$ be defined by
\begin{align}
  J_n(G,x) & := \int_{(\R^d)^n} \prod_{\{i,j\} \in E(G)}
\varphi(x_i-x_j) \ d(x_1, \dots, x_n).
\end{align}
By multiplying the integrand with a
$1=\varphi(x_i-x_j) + (1-\varphi(x_i-x_j))$  for each edge $\{i,j\}$
which is not contained in $E(G)$, we obtain that
\begin{align}
  J_n(G,x) & = \sum_{\substack{H \in \G_n \\ E(H) \supset E(G)}} I_n(H, x).
\end{align}
For example
\begin{align*}
  J_2\Bigg( \inclD{diagram_4_x1_x2_xy_12_n1y_n2y}\Bigg) =
  I_2\Bigg( \inclD{diagram_4_x1_x2_xy_12_n1y_n2y}\Bigg) +
  I_2\Bigg( \inclD{diagram_4_x1_x2_xy_12_1y_n2y}\Bigg) +
  I_2\Bigg( \inclD{diagram_4_x1_x2_xy_12_n1y_2y}\Bigg) +
  I_2\Bigg( \inclD{diagram_4_x1_x2_xy_12_1y_2y}\Bigg).
\end{align*}
By a M\"obius-inversion (see e.g.\ \cite{PeccTaqqu11}), we have
\begin{align}
  I_n(G,x) = \sum_{\substack{H \in \G_n \\ E(H) \supset E(G)}} (-1)^{|E(H)| - |E(G)|} J_n(H,x), \quad G \in \G_n.
\end{align}
This leads to the announced second representation for $p_n(x)$, namely
\begin{align*}
  p_n(x) & = \frac{1}{n!} \sum_{G \in \G_n} \pi_n(G) I_n(G,x) \\
  & = \frac{1}{n!} \sum_{G \in \G_n} \sum_{H \in \G_n} \ind\{E(H) \supset E(G)\}\pi_n(G)(-1)^{|E(H)| - |E(G)|} J_n(H,x) \\
  & = \frac{1}{n!} \sum_{H \in \G_n} J_n(H,x) \sum_{G \in \G_n}
  \ind\{E(G) \subset E(H)\}\pi_n(G)(-1)^{|E(H)| - |E(G)|} .
\end{align*}
In particular
\begin{align}\label{eq:kombinatorik:def_p_n_via_tilde_pi_n}
  p_n(x) & = \frac{1}{n!} \sum_{H \in \G_n} \kappa_n(H) J_n(H,x), \quad n \in \N
\end{align}
where
\begin{align}\label{eq:kombinatorik:def_tilde_pi_n}
  \kappa_n(H) := \sum_{G \in \G_n} \ind\{E(G) \subset E(H)\}\pi_n(G)(-1)^{|E(H)| - |E(G)|}.
\end{align}

A vertex $i\in [n]$ of a graph $G \in \G_n$ is called \emph{pivotal},
if any path from $0$ to $n+1$ contains $i$. The subset $\G_n^0 \subset \G_n$
of graphs which contain no pivotal vertex plays a significant
role for determining the coefficients $q_n(x)$ from $p_n(x)$ as the
next theorem shows.

\begin{theorem}\label{thm:repr_q_n}
The coefficients $q_n(x):=q_n(x,0)$ of the series
  representation of $Q_t$ at $t_0 = 0$ satisfy
  \begin{align*}
    q_n(x)=\frac{1}{n!}\sum_{H \in \G_n^0}\kappa_n(H)J_n(H, x),\quad x\in\R^d.
  \end{align*}
\end{theorem}

This means, that $q_n(x)$ differs from $p_n(x)$ only by the sum over
the graphs with pivotal vertices. The proof of Theorem
\ref{thm:repr_q_n} is based on some lemmas.


At first we define the concatenation of two graphs. For $n, m \in
\N_0$ let $G_1 \in \G_n$ and $G_2 \in \G_m$. The concatenation $G_1
\odot G_2 \in \G_{n+m+1}$ of $G_1$ and $G_2$ is constructed in the
following way:
\begin{enumerate}
  \item Relabel all nodes in $G_2$ with labels $\{n + 1, \dots, n + m + 2\}$ without changing the order.
  \item Define $V(G_1 \odot G_2) := V(G_1) \cup V(G_2)$.
  \item Define $E(G_1 \odot G_2) := E(G_1) \cup E(G_2)$.
\end{enumerate}
For example
\begin{align}
  \inclD{diagram_3_x1_nxy_1y} \odot \inclD{diagram_3_x1_xy_1y}
= \inclD{diagram_5_x1_nxz_nx2_nxy_1z_n12_n1y_z2_zy_2y}.
\end{align}
In other words, we only combine the end-vertex of $G_1$ and the
start-vertex of $G_2$ to a new vertex and adjust the labels.

\begin{lemma}\label{le:concatenation_multiplicative_pi}
  For $n, m \in \N_0$ and $G_1 \in \G_n$, $G_2 \in \G_m$ we have
  \begin{align}
    \pi_{n + m + 1}(G_1 \odot G_2) = \pi_n(G_1) \pi_m(G_2).
  \end{align}
\end{lemma}
\emph{Proof:} The vertex with label $n + 1$ is by construction
pivotal. If $0 \leftrightarrow n + m + 2$ in $(G_1 \odot G_2)|I$, then
$n + 1 \in I$. Hence
\begin{align}\notag
  \pi_{n+m+1}(G_1 \odot G_2) & = \sum_{I \subset [n + m + 1]} (-1)^{n + m + 1 - |I|}
\ind\{0\leftrightarrow n+m+2 \text{ in $(G_1 \odot G_2)|I$}\} \\ \notag
  & = \sum_{I_1 \subset [n]} \sum_{I_2 \subset [m]+n+1} (-1)^{n + m -
    |I_1| - |I_2|}
  \ind\{0 \leftrightarrow n+1 \text{ in $(G_1 \odot G_2)|I_1$}\} \\ \notag
  & \hspace{5cm}\times \ind\{n+1 \leftrightarrow n + m + 2 \text{ in $(G_1 \odot G_2)|I_2$}\} \\ \notag
  & = \sum_{I_1 \subset [n]} \sum_{I_2 \subset [m]}
(-1)^{n+m-|I_1|-|I_2|} \ind\{0 \leftrightarrow n + 1 \text{ in $G_1|I_1$}\} \\ \notag
  & \hspace{5cm} \times \ind\{0 \leftrightarrow m + 1 \text{ in $G_2|I_2$}\} \\
  \tag*{$\qed$}
  & = \pi_n(G_1) \pi_m(G_2).
\end{align}

\bigskip
\begin{lemma}\label{le:concatenation_multiplicative_kappa}
  For $n, m \in \N_0$ and $G_1 \in \G_n$, $G_2 \in \G_m$ we have
  \begin{align}
    \kappa_{n + m + 1}(G_1 \odot G_2) = \kappa_n(G_1) \kappa_m(G_2).
  \end{align}
\end{lemma}
\emph{Proof:} In every graph $H \in \G_{n + m + 1}$ with $E(H) \subset
E(G_1 \odot G_2)$ the vertex $n + 1$ is pivotal. Hence there are
uniquely determined Graphs $H_1 \in \G_n$ and $H_2 \in \G_m$ such that
$H = H_1 \odot H_2$. The graph $H_1$ consists of $0$, $n+1$ and all
vertices lying "in front" of $n+1$, whereas $H_2$ is a relabeled
version of the subgraph of $H$ which consists of $n+1$, $n + m + 2$
and all vertices lying "behind" $n+1$. Hence by Lemma
\ref{le:concatenation_multiplicative_pi}
\begin{align}\notag
  \kappa_{n + m + 1}&(G_1 \odot G_2) \\ \notag
  & = \sum_{H \in \G_{n + m + 1}} \ind\{E(H)
\subset E(G_1 \odot G_2)\} \pi_{n + m + 1}(H) (-1)^{|E(H)| - |E(G_1 \odot G_2)|} \\ \notag
  & = \sum_{H_1 \in \G_n} \sum_{H_2 \in \G_m}\ind\{E(H_1\odot H_2)
\subset E(G_1 \odot G_2)\}\pi_{n+m+1}(H_1\odot H_2) \\ \notag
  & \hspace{3cm} \times (-1)^{|E(H_1 \odot H_2)| - |E(G_1 \odot G_2)|} \\ \notag
  & = \sum_{H_1 \in \G_n} \sum_{H_2 \in \G_m} \ind\{E(H_1)
\subset E(G_1)\} \ind\{E(H_2) \subset E(G_2)\} \\ \notag
  & \hspace{3cm} \times \pi_{n}(H_1) \pi_m(H_2) (-1)^{|E(H_1)| + |E(H_2)| - |E(G_1)| - |E(G_2)|} \\
\tag*{$\qed$}
  & = \kappa_n(G_1) \kappa_m(G_2).
\end{align}

\bigskip

\begin{lemma}\label{le:concatenation_convolution_J}
  For $n, m \in \N_0$ and $G_1 \in \G_n$, $G_2 \in \G_m$ we have
  \begin{align}
    J_n(G_1) \ast J_m(G_2) = J_{n + m + 1}(G_1 \odot G_2).
  \end{align}
\end{lemma}
\emph{Proof:} For all $x \in \R^d$, we have
\begin{align*}
  (J_n(G_1) \ast J_m(G_2))(x) & = \int_{\R^d} \int_{(\R^d)^n} \prod_{\{i,j\} \in E(G_1)}
\varphi(x_i-x_j) \ d(x_1, \dots, x_n) \\
  & \hspace{.6cm} \times \int_{(\R^d)^m} \prod_{\{i,j\} \in E(G_2)}
  \varphi(y_i-y_j)\ d(y_1, \dots, y_m) \ dx_{n+1},
\end{align*}
where $x_0:=0$ and $y_0 := 0$,
$y_{m+1}:=x-x_{n+1}$. By translation invariance, nothing changes, if we
redefine $y_0 := x_{n+1}$ and $y_{m+1} := x$. If we apply Fubinis Theorem
and rename the integration variables in the same way as we renamed the vertex labels
in the definition of "concatenation" we obtain
\begin{align*}
  & \quad \ (J_n(G_1) \ast J_m(G_2))(x) \\
  & = \int_{\R^d} \int_{(\R^d)^n} \int_{(\R^d)^m}\prod_{\{i,j\} \in E(G_1 \odot G_2)}
\varphi(x_i-x_j) \  d(x_1, \dots, x_n)\ d(x_{n+2}, \dots, x_{n+m+1}) \ dx_{n+1} \\
  & = J_{n+m+1} (G_1 \odot G_2)(x).
\end{align*}
\qed

\bigskip

We are now in the position to prove Theorem \ref{thm:repr_q_n}.

\bigskip \emph{Proof of Theorem \ref{thm:repr_q_n}:} For $t_0 = 0$
equation \eqref{eq:analyticity_of_Q_t:rec1} and
\eqref{eq:analyticity_of_Q_t:rec2} simplify to $q_0 = p_0$ and
\begin{align}\label{eq:thm_repr_q_n:1}
  q_{n+1} = p_{n+1} - \sum_{k = 0}^n q_{n-k} \ast p_{k}, \quad n \in \N_0.
\end{align}
We will use this for an induction over $n$. First we observe, that
trivially $\G_0 = G_0^0$, as the graph $G_0$ that connects $0$ and $1$
with a single bond is the only element of $\G_0$. Hence
\begin{align}
  \frac{1}{0!} \sum_{G \in \G_0^0} \kappa_0(G) J_0(G, x) = J_0(G_0, x) = p_0 = q_0.
\end{align}
For the induction step we define $\G_n^{>0} := \G_n \setminus
\G_n^0$. If follows from
\eqref{eq:kombinatorik:def_p_n_via_tilde_pi_n} that
\begin{align}
  p_{n+1} = \frac{1}{(n+1)!} \sum_{G \in \G_{n+1}^0} \kappa_{n+1}(G)
  J_{n+1}(G) + \frac{1}{(n+1)!} \sum_{G \in \G_{n+1}^{>0}}
  \kappa_{n+1}(G) J_{n+1}(G).
\end{align}
Hence by \eqref{eq:thm_repr_q_n:1} it is enough to show
\begin{align}\label{eq:thm_repr_q_n:2}
  \frac{1}{(n+1)!} \sum_{G \in \G_{n+1}^{>0}} \kappa_{n+1}(G) J_{n+1}(G) = \sum_{k = 0}^n p_k \ast q_{n-k}.
\end{align}
We use the induction hypothesis,
\eqref{eq:kombinatorik:def_p_n_via_tilde_pi_n} and the Lemmas
\ref{le:concatenation_multiplicative_kappa} and
\ref{le:concatenation_convolution_J} to obtain
\begin{align}
 \sum_{k=0}^n& p_k \ast q_{n-k} \nonumber \\
  & = \sum_{k = 0}^n \left( \frac{1}{k!} \sum_{G_1 \in \G_k}
    \kappa_k(G_1) J_k(G_1) \right)
  \ast \bigg( \frac{1}{(n-k)!} \sum_{G_2 \in \G_{n-k}^0} \kappa_{n-k}(G_2) J_{n-k}(G_2) \bigg) \nonumber \\
  & = \sum_{k = 0}^n \frac{1}{k!(n-k)!} \sum_{G_1 \in \G_k} \sum_{G_2
    \in \G_{n-k}^0} \kappa_{n+1}(G_1 \odot G_2) J_{n+1}(G_1 \odot
  G_2) \label{eq:kombinatorische_darstellung_der_q_n:1}.
\end{align}
Finally we have a look at the left-hand side of
\eqref{eq:thm_repr_q_n:2}. Let $G \in \G_{n+1}^0$. Each path from $0$
to $n+2$ runs through the pivotal vertices of $G$ in the same
order. Let $v \in [n+1]$ be the last of these pivotal vertices. We
define the set of graphs $\H_k \subset \G_{n+1}^0$, $k \in \{0, \dots,
n\}$ with the following properties:
\begin{enumerate}
  \item Each $H \in \H_k$ contains at least one pivotal vertex.
  \item The vertex with label $k + 1$ is the last pivotal vertex in each $H \in \H_k$.
  \item The $k$ vertices $\{1, \dots, k\}$ lie in front of the vertex $k + 1$.
  \item The $n - k$ vertices $\{k + 2, \dots, n+1\}$ lie behind the vertex $k + 1$.
\end{enumerate}
We partition the set $[n+1]$ of vertices in each $G \in \G_{n+1}^0$
into three sets $M_1$, $M_2$ and $M_3$. The set $M_1$ contains all
vertices that lie in front of the last pivotal vertex of $G$. The set
$M_2$ contains only the last pivotal vertex and $M_3$ contains the
remaining vertices. Now we relabel the vertices in $G$ to obtain a
graph $\tilde G$ in the following way: The vertices in $M_1$ are
labeled with the numbers $1, \dots, |M_1|$ without changing the
order. The vertex in $M_2$ is labeled $|M_1| + 1$ and the vertices in
$M_3$ are labeled with the numbers $|M_1| + 2, \dots, n + 1$, again
without changing the order. By construction we have $\tilde G \in
\H_{|M|_1}$ but $\kappa(G) = \kappa(\tilde G)$ and $J_{n+1}(G) =
J_{n+1}(\tilde G)$. There are exactly $\binom{n+1}{k, n-k, 1}$ Graphs
in $\G_{n+1}^0$, which become the same $\tilde G$ by this relabeling
procedure. Hence we have for the left-hand side of
\eqref{eq:thm_repr_q_n:2}
\begin{align*}
  \frac{1}{(n+1)!} \sum_{G \in \G_{n+1}^{>0}} \kappa_{n+1}(G) J_{n+1}(G)
& = \frac{1}{(n+1)!} \sum_{k = 0}^n \sum_{H \in \H_k} \binom{n+1}{k, n-k, 1} \kappa_{n+1}(H) J_{n+1}(H) \\
  & = \sum_{k = 0}^n \sum_{H \in \H_k} \frac{1}{k! (n-k)!}
  \kappa_{n+1}(H) J_{n+1}(H)
\end{align*}
which is equal to \eqref{eq:kombinatorische_darstellung_der_q_n:1} due to the definition of the concatenation.\qed

\bigskip

\section{Appendix: Palm distributions}
In this paper all random elements are defined on a measurable space
$(\Omega,\mathcal{A})$ equipped with a {\em measurable flow}
$\theta_x\colon\Omega \to \Omega$, $x\in \R^d$. This is a family
of measurable mappings such that $(\omega ,x)\mapsto \theta_x
\omega$ is measurable, $\theta_0$ is the identity on $\Omega$ and
\begin{align}\label{flow}
\theta_x \circ \theta_y =\theta_{x+y},\quad x,y\in \R^d,
\end{align}
where $\circ$ denotes composition.
We may think of $\theta_x\omega$ as of
$\omega$ {\em shifted} by the vector $-x$.
We fix a probability measure
$\P$ on $(\Omega,\mathcal{A})$ and assume that it is
{\em stationary}, that is
$$
\P\circ\theta_x=\P,\quad x\in\R^d,
$$
where $\theta_x$ is interpreted as a mapping from $\mathcal{A}$ to $\mathcal{A}$
in the usual way:
$$
\theta_xA:=\{\theta_x\omega:\omega\in A\},\quad A\in\mathcal{A},\, x\in\R^d.
$$

Let  $\bN(\R^d)$ denote the space of all
locally finite subsets $\mu$ of $\R^d$. Hence $\mu\in\bN(\R^d)$
iff $\mu\cap B$ is finite for each bounded set.
For each $\mu\in \bN(\R^d)$ and each $B\subset\R^d$ we write
$\mu(B):=|\mu\cap B|$ for the number of points of $\mu$ lying in $B$.
As usual we equip $\bN(\R^d)$ with the smallest $\sigma$-field $\mathcal{N}$
making the mappings $\mu\mapsto\mu(B)$ measurable for all
$B$ in the Borel $\sigma$-field $\cB^d$ on $\R^d$.


A {\em point process} on $\R^d$ is
a measurable mapping $\eta\colon\Omega\rightarrow\mathbf{N}(\R^d)$.
It is called {\em invariant} (or {\em stationary})
\begin{align}\label{invariant}
\eta(\omega,B+x)=\eta(\theta_x\omega,B),\quad \omega\in \Omega,x\in \R^d,
B\in\mathcal{B}^d.
\end{align}
Let $\eta$ be an invariant point process.
The {\em intensity} of $\eta$ is the
number $\gamma_\eta:=\E\eta([0,1]^d)$. If the latter
is positive and finite we can define
the probability measure
\begin{align} \label{Palm}
\P^0_\eta(A):=\gamma_\eta^{-1}\int\sum_{x\in\eta(\omega)} \I\{\theta_x\omega\in A, x\in[0,1]^d\}\,
\P(d\omega), \quad A\in\mathcal{A}.
\end{align}
This {\em Palm probability measure} of $\eta$
satisfies the {\em refined Campbell formula}
\begin{align}\label{refC}
\int\sum_{x\in\eta(\omega)} f(\theta_x\omega,x)\, \P(d\omega)=
\gamma_\eta\iint f(\omega,x)\,dx\,\P^0_\eta(d\omega)
\end{align}
for all measurable $f\colon\Omega\times\R^d\to [0,\infty)$, where
$dx$ refers to integration with respect to Lebesgue measure on $\R^d$.
Using standard conventions we write this as
\begin{align} \label{21}
\E\sum_{x\in\eta} f(\theta_x,x)=\gamma_\eta\,\E^0_\eta\int f(\theta_0,x)\,dx,
\end{align}
where $\E^0_\eta$ denote integration with respect to $\P^0_\eta$.
The measure $\P^0_\eta$ is concentrated on the measurable set
$\Omega_0$ of all $\omega\in\Omega$ such that the origin $0$ is in $\eta(\omega)$.
The {\em Palm distribution} of $\eta$ is the distribution $\P^0_\eta(\eta\in\cdot)$
of $\eta$ under $\P^0_\eta$. It is concentrated on the measurable set
of all $\mu\in\mathbf{N}(\R^d)$ such that $0\in\mu$.
The number $\P_\eta^0(A)$ can be interpreted as the conditional probability of
$A\in\mathcal{A}$ given that $\eta$ has a point at the origin.

The following result (Neveu's exchange formula) is a versatile
tool of Palm theory.

\begin{proposition}\label{tneveu} Let $\eta,\eta'$ be two invariant point
processes with finite intensities and
let $f\colon\Omega\times\R^d\to[0,\infty)$ be measurable. Then
\begin{align}\label{A2}
  \gamma_\eta\,\E^0_\eta \sum_{x\in\eta'} f(\theta_0,x)
=\gamma_{\eta'}\,\E^0_{\eta'} \sum_{x\in\eta} f(\theta_x,-x).
\end{align}
This remains true for any measurable $f\colon\Omega\times\R^d\to \R$ with
$\E^0_{\eta'} \int |f(\theta_x,-x)|\,\eta(dx)<\infty$.
\end{proposition}

Let $\eta$ be a point process on $\R^d$ and $n\in\N$.
The $n$-th {\em factorial moment measure}
$\alpha^{(n)}$ of $\eta$ is the measure on $(\R^d)^n$, defined by
\begin{align}\label{factorial}
\alpha^{(n)}:=\E\sumneq_{x_1,\ldots,x_n\in\eta}\I\{(x_1,\ldots,x_n)\in\cdot\},
\end{align}
where $\sum^{\ne}$ denotes summation over all ordered $n$-tuples of distinct elements
of $\eta$.
Assume now  that $\eta$ is an invariant point process with a positive and
finite intensity. Assume also
that $\alpha^{(2)}$ is locally finite (finite on bounded
Borel sets) and absolutely continuous with respect to
Lebesgue measure, that is (using also stationarity)
\begin{align}\label{paircorr}
\alpha^{(2)}(d(x,y))=\gamma_\eta^2\rho(y-x)\,d(x,y)
\end{align}
for a locally integrable measurable $\rho\colon\R^d\to\R$. The latter is
the {\em pair correlation function} of $\eta$.
The {\em two-point Palm probability measures} of
$\eta$ is a family $\{\P^{x,y}_\eta:x,y\in\R^d\}$ of probability
measures on $(\Omega,\mathcal{A})$ such that
$(x,y)\mapsto  \P^{x,y}_\eta(A)$ is measurable for all $A\in\mathcal{A}$
and
\begin{align}\label{A2Palm}
\E\sumneq_{x,y\in\eta} f(\theta_0,x,y)
=\gamma_\eta^2\int \E^{x,y}_\eta f(\theta_0,x,y)\rho(y-x)\,d(x,y)
\end{align}
for all measurable $f\colon \Omega\times\R^d\times\R^d\to[0,\infty)$,
where $\E^{x,y}_\eta$ denotes expectation with respect to $\P^{x,y}_\eta$.
These Palm distributions can be constructed as follows.
An easy calculation shows that for any $B\in\mathcal{B}^d$
\begin{align*}
  \E^0_\eta\sum_{x\in\eta\setminus\{0\}} \I\{x\in B\}=\gamma_\eta^2\int\I\{x\in B\}\rho(x)\,dx.
\end{align*}
We now assume that $(\Omega,\mathcal{A})$ is a {\em Borel space}.
This very weak assumption can be made without restricting generality.
By a standard disintegration technique we can then
find a family $\{\P^{0,x}_\eta:x\in\R^d\}$ of probability measures on $(\Omega,\mathcal{A})$
such that
\begin{align}\label{A5}
\E^0_\eta\int f(\theta_0,x)\,\eta(dx)
=\E^0_\eta f(\theta_0,0)+\gamma_\eta\int\E^{0,x}_\eta f(\theta_0,x)\rho(x)\,dx
\end{align}
for all measurable $f\colon \Omega\times\R^d\to[0,\infty)$.
We can then define
\begin{align*}
\P^{x,y}_\eta(A):=\P^{0,y-x}(\theta_xA),\quad x,y\in\R^d,\,A\in\mathcal{A}.
\end{align*}
This means that
\begin{align}\label{A3}
  \E^{x,y}_\eta F= \E^{0,y-x}_\eta F\circ\theta_{-x},\quad x,y\in\R^d,
\end{align}
for all measurable $F\colon\Omega\to[0,\infty)$. Using the refined
Campbell theorem \eqref{refC} and \eqref{A5} it is then not hard
to check that \eqref{A2Palm} holds.
It is also easy to see that $\P^{x,y}_\eta(x,y\in\eta)=1$
for $\alpha^{(2)}$-a.e.\ $(x,y)\in \R^d\times\R^d$.
The number $\P^{x,y}_\eta(A)$ can be interpreted
as  probability of $A\in\mathcal{A}$ given that $\eta$ has points
at $x$ and $y$.

Let us now assume that $\eta\equiv \eta_t$ is a stationary
Poisson process of intensity $t> 0$.
The multivariate Mecke equation (see e.g.\ \cite[Corollary 3.2.3]{SW08}
or \cite[Theorem 4.4]{LastPenrose16}) states for any $n\in\N$ and any measurable function
$f\colon\bN(\R^d) \times (\R^d)^n \to [0,\infty )$ that
\begin{align}\label{Mecke}
\E\sumneq_{x_1,\ldots,x_n\in\eta} f(\eta_t,x_1,\dots,x_n)
=t^n\int\E f\big(\eta_t\cup\{x_1,\ldots,x_n\},x_1,\dots,x_n\big)\,d(x_1,\dots, x_n).
\end{align}
The case $n=1$ easily implies (together with stationarity of $\eta_t$)
that the Palm distribution of $\eta_t$ is given by
\begin{align*}
\P^0_{\eta_t}(\eta_t\in\cdot)=\P(\eta_t\cup\{0\}\in\cdot).
\end{align*}
For $n=2$ we obtain from \eqref{Mecke} that the pair correlation function $\rho_t$ of $\eta_t$
satisfies $\rho_t\equiv 1$ and that, moreover,
\begin{align*}
\P^{x,y}_{\eta_t}(\eta_t\in\cdot)=\P(\eta_t\cup\{x,y\}\in\cdot)
\end{align*}
for almost every $(x,y)$ with respect to Lebesgue measure on $(\R^d)^2$.

In this paper we work also with point processes on a metric space
$\X$ different from $\R^d$. These are a random elements of the space
$\bN(\X)$ of all
integer-valued locally finite measures $\mu$
on $\X$ equipped with
the smallest $\sigma$-field
making the mappings $\mu\mapsto\mu(B)$ measurable for all
$B$ in the Borel $\sigma$-field on $\X$.
For more details on point processes we refer to \cite{LastPenrose16,SW08}.
A survey of Palm theory can be found in \cite{Last10}.

\bigskip
\noindent {\bf Acknowledgments:} We wish to thank Salvatore Torquato for
drawing our attention to the Ornstein-Zernike equation and its relevance
for percolation theory. This work was supported
by the German Research Foundation (DFG) through the research unit 
``Geometry and Physics of Spatial Random System'' 
under the grant LA 965/7-2.

\end{document}